\documentclass[hidelinks,onefignum,onetabnum,dvipsnames]{siamart220329}

\def\ArticleTitle{A multigrid reduction framework for domains with symmetries}
\def\ShortArticleTitle{An AMG reduction framework for symmetric domains}

\usepackage{amsfonts}
\usepackage{graphicx}
\usepackage{epstopdf}
\usepackage{array,multirow}

\usepackage{tabularx,booktabs}
\newcolumntype{Y}{>{\centering\arraybackslash}X}
\newcolumntype{C}[1]{>{\centering\arraybackslash}m{#1}}
\usepackage{tabularx}

\ifpdf
  \DeclareGraphicsExtensions{.eps,.pdf,.png,.jpg}
\else
  \DeclareGraphicsExtensions{.eps}
\fi

\newsiamremark{remark}{Remark}
\newsiamremark{hypothesis}{Hypothesis}
\crefname{hypothesis}{Hypothesis}{Hypotheses}
\newsiamthm{claim}{Claim}

\headers{\ShortArticleTitle}{A.~Alsalti-Baldellou, C.~Janna, X.~\'Alvarez-Farr\'e, and F.~X.~Trias}

\title{\ArticleTitle\thanks{Accepted for publication.
\funding{A.A.B., X.A.F. and F.X.T. were supported by the SIMEX project (PID2022-142174OB-I00) of \emph{Ministerio de Ciencia e Innovaci\'{o}n} and the RETOtwin project (PDC2021-120970-I00) of \emph{Ministerio de Econom\'{i}a y Competitividad}, Spain.
A.A.B. and C.J. were supported by the project ``National Centre for HPC, Big Data and Quantum Computing'', CN00000013 (approvato nell'ambito del Bando M42C – Investimento 1.4 – Avvisto ``Centri Nazionali'' – D.D.~n.~3138 del 16.12.2021, ammesso a finanziamento con Decreto del MUR n.~1031 del 17.06.2022).}}}

\author{
\`Adel Alsalti-Baldellou\footnotemark[1] \footnotemark[2]
\and
Carlo Janna\footnotemark[3] \footnotemark[4]
\and
Xavier \'Alvarez-Farr\'e\footnotemark[5]
\and
F. Xavier Trias\footnotemark[2]
}

\usepackage{amsopn}

\usepackage{xcolor}

\usepackage{etoolbox}

\newcommand{\ie}{\textit{i.e.},}
\newcommand{\eg}{\textit{e.g.},}

\newbool{xbcfd}
\setbool{xbcfd}{true}
\newcommand{\xcfd}{\ifbool{xbcfd}{computational fluid dynamics (CFD)\global\setbool{xbcfd}{false}}{CFD}}

\newbool{xbdns}
\setbool{xbdns}{true}
\newcommand{\xdns}{\ifbool{xbdns}{direct numerical simulation (DNS)\global\setbool{xbdns}{false}}{DNS}}

\newbool{xbles}
\setbool{xbles}{true}
\newcommand{\xles}{\ifbool{xbles}{large-eddy simulation (LES)\global\setbool{xbles}{false}}{LES}}

\newbool{xbspmm}
\setbool{xbspmm}{true}
\newcommand{\xspmm}{\ifbool{xbspmm}{sparse matrix-matrix product (\texttt{SpMM})\global\setbool{xbspmm}{false}}{\texttt{SpMM}}}

\newbool{xbspmv}
\setbool{xbspmv}{true}
\newcommand{\xspmv}{\ifbool{xbspmv}{sparse matrix-vector product (\texttt{SpMV})\global\setbool{xbspmv}{false}}{\texttt{SpMV}}}

\newbool{xbpcg}
\setbool{xbpcg}{true}
\newcommand{\xpcg}{\ifbool{xbpcg}{Preconditioned Conjugate Gradient (PCG)\global\setbool{xbpcg}{false}}{PCG}}

\newbool{xbamg}
\setbool{xbamg}{true}
\newcommand{\xamg}{\ifbool{xbamg}{Algebraic Multigrid (AMG)\global\setbool{xbamg}{false}}{AMG}}

\newbool{xblrcamg}
\setbool{xblrcamg}{true}
\newcommand{\xlrcamg}{\ifbool{xblrcamg}{Low-Rank Corrected \xamg~(LRCAMG($k$))\global\setbool{xblrcamg}{false}}{LRCAMG($k$)}}

\newbool{xbfsai}
\setbool{xbfsai}{true}
\newcommand{\xfsai}{\ifbool{xbfsai}{Factored Sparse Approximate Inverse (FSAI)\global\setbool{xbfsai}{false}}{FSAI}}

\newbool{xblrcfsai}
\setbool{xblrcfsai}{true}
\newcommand{\xlrcfsai}{\ifbool{xblrcfsai}{Low-Rank Corrected \xfsai~(LRCFSAI($k$))\global\setbool{xblrcfsai}{false}}{LRCFSAI($k$)}}

\def\K{K}
\def\bK{\bar{\K}}
\def\B{B}
\def\bB{\bar{\B}}
\def\Bt{B^T}
\def\bBt{\bar{\B}^T}
\def\C{C}
\def\bC{\bar{\C}}
\def\S{S}
\def\hS{\hat{\S}}
\def\bS{\bar{\S}}
\def\ninn{n_\text{inn}}
\def\nifc{n_\text{ifc}}

\DeclareMathOperator{\rank}{rank}
\DeclareMathOperator{\range}{range}

\def\R{\mathbb{R}}
\def\Rn{\mathbb{R}^n}
\def\Rnn{\mathbb{R}^{n \times n}}
\def\Id{\mathbb{I}}
\def\H{H}
\def\Hinn{\H_\text{inn}}
\def\Hout{\H_\text{out}}
\def\A{A}
\def\hA{\hat{\A}}
\def\hAi{\hat{\A}_i}
\def\Ainn{\A_\text{inn}}

\def\Aouti{\A_{\text{out},i}}
\def\G{G}

\def\Ginn{\G_\text{inn}}
\def\GK{\G_{\K}}
\def\GbK{\G_{\bK}}
\def\P{Q}
\def\Ps{\P_s}
\def\H{H}

\def\AMGinn{M_\text{inn}^{-1}}
\def\AMGA{M_{\A}^{-1}}
\def\AMGK{M_{\K}^{-1}}
\def\AMGbK{M_{\bK}^{-1}}
\def\AMGS{M_\text{AMGS}^{-1}}
\newcommand{\dotcup}{\ensuremath{\mathaccent\cdot\cup}}

\def\GbC{\G_{\bC}}

\def\W{W}
\def\bW{\bar{\W}}

\usepackage{algorithm}
\usepackage{algorithmicx}
\usepackage{algpseudocode}
\usepackage{xfrac}
\usepackage{mathtools}
\usepackage{pgfplotstable}
\usepackage[normalem]{ulem}
\usepackage{hyperref}
\usepackage{subcaption}
\usepackage{xifthen}

\ifpdf
\hypersetup{
  pdftitle={\ArticleTitle},
  pdfauthor={A.~Alsalti-Baldellou, C.~Janna, X.~\'Alvarez-Farr\'e, and F.~X.~Trias}
}
\fi

\begin{document}

\renewcommand{\thefootnote}{\fnsymbol{footnote}}
\footnotetext[1] {Department of Information Engineering,
University of Padova, Via Giovanni Gradenigo, 6b, 35131 Padova PD, Italy (\email{adel.alsaltibaldellou@unipd.it}).}
\footnotetext[2] {Heat and Mass Transfer Technological Center, Universitat Polit\`{e}cnica de Catalunya -- BarcelonaTech, Carrer de Colom, 11, 08222 Terrassa, Spain (\email{francesc.xavier.trias@upc.edu}).}
\footnotetext[3] {Department of Civil, Environmental and Architectural Engineering, University of Padova, Via Francesco Marzolo, 9, 35131 Padova PD, Italy (\email{carlo.janna@unipd.it}).}
\footnotetext[4] {M$^3$E S.r.l., Via Giambellino, 7, 35129 Padova PD, Italy.}
\footnotetext[5] {High-Performance Computing and Visualization Team, SURF, Science Park 140, 1098 XG Amsterdam, The Netherlands (\email{xavier.alvarezfarre@surf.nl}).}
\renewcommand{\thefootnote}{\arabic{footnote}}

\maketitle

\begin{abstract}
Divergence constraints are present in the governing equations of numerous physical phenomena, and they usually lead to a Poisson equation whose solution represents a bottleneck in many simulation codes.
Algebraic Multigrid (AMG) is arguably the most powerful preconditioner for Poisson's equation, and its effectiveness results from the complementary roles played by the smoother, responsible for damping high-frequency error components, and the coarse-grid correction, which in turn reduces low-frequency modes.
This work presents several strategies to make AMG more compute-intensive by leveraging reflection, translational and rotational symmetries.
AMGR, our final proposal, does not require boundary conditions to be symmetric, therefore applying to a broad range of academic and industrial configurations.
It is based on a multigrid reduction framework that introduces an aggressive coarsening to the multigrid hierarchy, reducing the memory footprint, setup and application costs of the top-level smoother.
While preserving AMG's excellent convergence, AMGR allows replacing the standard sparse matrix-vector product with the more compute-intensive sparse matrix-matrix product, yielding significant accelerations.
Numerical experiments on industrial CFD applications demonstrated up to 70\% speed-ups when solving Poisson's equation with AMGR instead of AMG.
Additionally, strong and weak scalability analyses revealed no significant degradation.

\end{abstract}



\section{Introduction}\label{sec:intro}

Divergence constraints are prevalent in physical problems, often following fundamental conservation principles like mass or electrical charge conservation.
Such constraints lead to a Poisson equation that plays a fundamental role in many areas of science and engineering, such as \xcfd, linear elasticity, and electrostatics.
Indeed, the solution of the associated linear system is usually the most computationally intensive part of scientific simulation codes, and the design and implementation of Poisson solvers are far from straightforward due to the interplay of numerical and hardware challenges.

Historically, the most efficient way to solve Poisson's problems is through iterative methods based on Krylov subspaces \cite{Axe94,Saa03,VdV03}, whose implementation is simple and easily parallelisable, requiring only basic linear algebra operations.
Namely, matrix by vector products, scalar products and vector updates.
However, iterative linear solvers must be properly preconditioned to be effective \cite{SaaVor00,Ben02}.
The choice, design and implementation of such preconditioners are not trivial, and it is one of the most active research fields in numerical analysis.
Preconditioners based on incomplete factorisations were very popular in the early days of numerical linear algebra.
The first papers suggesting this approach were \cite{MeiVdV77,Ker78}, where the factorisation takes place with no fill-in, \ie~the pattern of the Cholesky factor, $L$, equals the lower pattern of $A$, with $A \simeq LL^T$.
Subsequently, more advanced and effective alternatives with dynamic fill-in control were proposed; see, for instance, \cite{Saa94,LinMor99}.

However, the sequential nature of such methods and the increasing availability of parallel computers in the late '90s made them lose ground against alternatives with higher degrees of parallelism.
For instance, the application of preconditioners based on approximate inverses solely relies on the \xspmv, an easily parallelisable operation.
In some cases, even their construction is reasonably concurrent.
The most prominent variants are AINV~\cite{BenMey95,BenMeyTum96,BenCulTum00},
SPAI~\cite{GroHuc97} and the \xfsai~\cite{KolYer93,JanFerSarGam15}.
While approximate inverses provide a high degree of parallelism, they are not {\em optimal} in the sense that when augmenting the mesh resolution (hence increasing the linear system size), the problem becomes more ill-conditioned and more iterations are required to reach the same accuracy.

This problem worsens nowadays, as extreme-scale linear systems must be solved on massively parallel supercomputers, and single-level preconditioners such as incomplete factorisations or approximate inverses generally require excessive iterations.
The problem of scalability is overcome with multilevel preconditioners like Geometric Multigrid (GMG) or \xamg~\cite{Tro00,Vas08,XuZik17}.
Thanks to the interplay between smoother and coarse-grid correction, these methods often solve a given PDE with a number of iterations independent of the mesh size.
Many freely available multigrid packages exist.
To cite a few, Hypre \cite{hypreURL}, Trilinos \cite{mueluURL}, and PETSc \cite{petscURL} provide very effective implementations with excellent scalability.

To develop efficient and scalable solvers, it is necessary to identify the limitations of current computing devices and develop algorithms that overcome them.
For instance, the low arithmetic intensity of most sparse linear algebra kernels motivated strategies like using mixed precision \cite{Baboulin2009} or applying more compute-intensive algorithms \cite{Ibeid2020}.
Similarly, the large memory to network bandwidth ratio led to implementations that hide or completely avoid inter-node communications \cite{Laut2022,Nayak2021,Sanan2016}.
Additionally, the limited available memory resulted in approaches like exploiting data sparsity \cite{Ghysels2016,Amestoy2015}.

This work presents multiple strategies leading to a multigrid reduction framework that mitigates several of the aforementioned computational challenges.
Namely, we show that given an arbitrarily complex geometry presenting reflection, translational or rotational symmetries, it is possible to apply a consistent ordering that makes the coefficient matrix (and preconditioners) conform to regular block structures regardless of the boundary conditions.
This, in turn, enables replacing the standard \xspmv{} with the more compute-intensive \xspmm.

Similar strategies were used in the context of \xcfd{} simulations~\cite{Shishkina2009}, but they were limited to reflection symmetries, and their main goal was to reduce the memory footprint~\cite{Alsalti2023b,Gorobets2011}.
In \cite{Alsalti2023a}, such strategies were generalised by leveraging an arbitrary number of reflection symmetries to block diagonalise the discrete Laplacian, which improved iterative methods' convergence and granted several computational advantages.
Nevertheless, the resulting subsystems were not identical, preventing the use of \xspmm{} on complex preconditioners.
In \cite{Alsalti2023c}, low-rank corrections were introduced to make \xfsai{} compatible with \xspmm{}, resulting in a much faster and significantly lighter variant.
However, given the poor weak scaling of single-level preconditioners, the present work extends past strategies to multilevel preconditioning.
Namely, to \xamg, arguably the most powerful method for preconditioning Poisson's equation.

Extending \xspmm's use to \xamg{} is far from direct, and for completeness, intermediate unsatisfactory approaches leading to AMGR, our final proposal, are included in this work.
AMGR relies on a multigrid reduction framework that introduces an aggressive coarsening to the top of the multigrid hierarchy, reducing the memory footprint, setup and application costs of the top-level smoother, which replaces \xspmv{} with \xspmm.
Subsequent levels are purely algebraic.
Apart from reflections, AMGR also leverages translational and rotational symmetries and, most notably, is compatible with problems with asymmetric boundary conditions, therefore applying to a significantly broader range of industrial applications, on which it proved up to 70\% faster than the standard \xamg{} thanks to preserving its excellent convergence while enjoying the computational advantages of \xspmm.

Without loss of generality, the targeted applications are incompressible \xcfd~simulations.
In particular, \xdns~and \xles~of turbulent flows, which are essential for many areas of engineering, such as energy production and environmental monitoring.
In such areas, spatial symmetries are commonly exhibited.
For instance, most vehicles have a central reflection symmetry~\cite{Aljure2018,Choi2014}.
Cases with two reflection symmetries are also usual, including jets~\cite{Duponcheel2021} and building simulations~\cite{Morozova2020}.
Further gains can be attained by exploiting translational or rotational symmetries.
They are present in turbomachinery~\cite{Wang2023}, heat exchangers~\cite{Lotfi2014,Paniagua2014}, nuclear reactors~\cite{Fang2021,Merzari2020}, or recurring structures, such as arrays of buildings~\cite{Hang2012} or wind farms~\cite{Calaf2010}.

The remaining sections are organised as follows.
\Cref{sec:sympois} recalls a strategy for exploiting reflection symmetries through low-rank corrections, and extends it to \xamg. Then, \cref{sec:innifc} introduces an ordering able to cope with asymmetric boundary conditions, as well as translational and rotational symmetries, and ultimately leading to AMGR.
\Cref{sec:impl} discusses the parallel implementation, \cref{sec:experiments} presents numerical experiments on industrial \xcfd{} applications, and \cref{sec:concl} gives some concluding remarks.

\section{Symmetry-aware spatial discretisation}\label{sec:sympois}

\begin{figure}[tbp]
  \centering
  \includegraphics[trim = 5.8cm 20.6cm 8.2cm 8.1cm,clip,width=\linewidth]{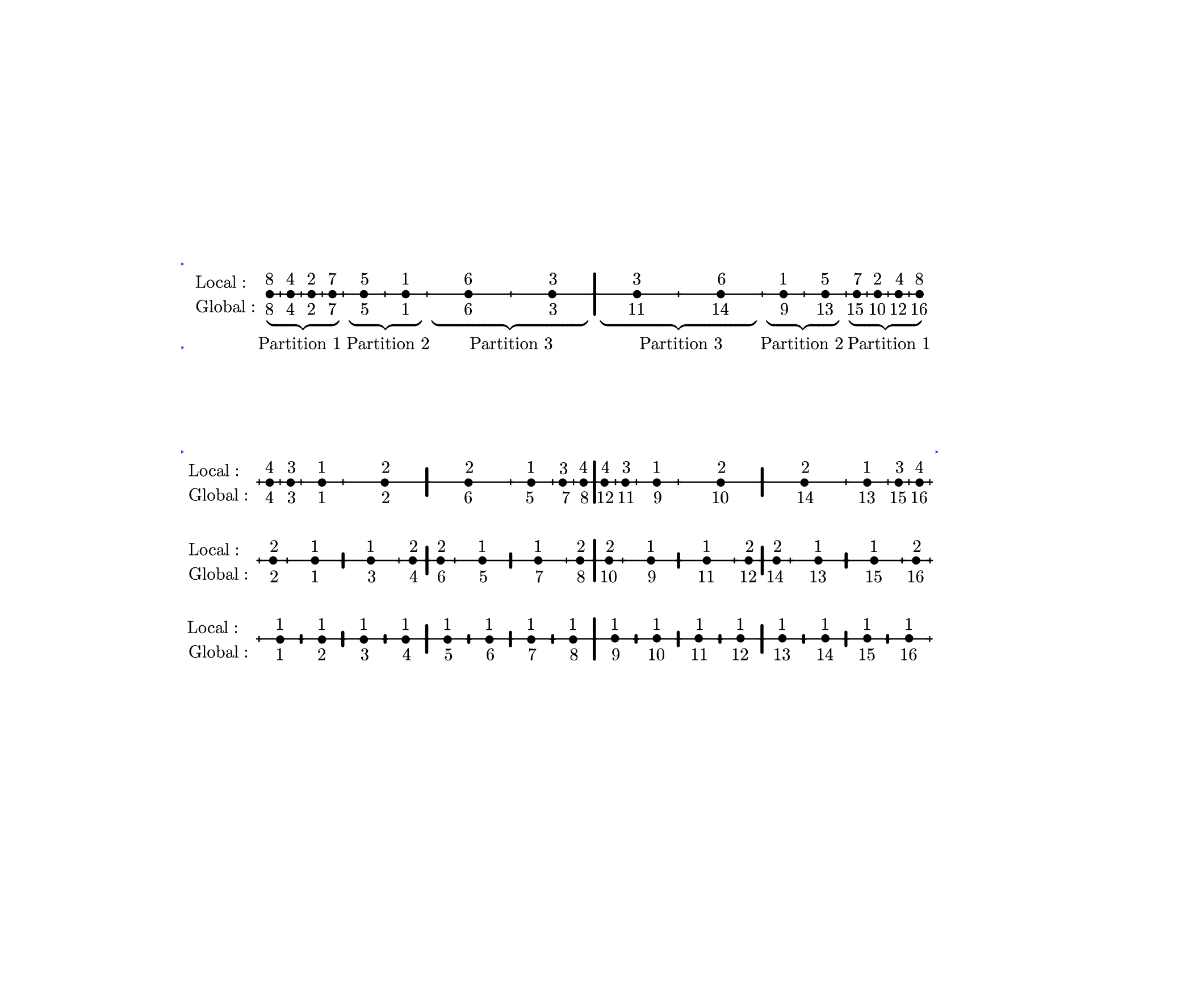}
  \caption{Single-symmetry 1D mesh with a random mirrored ordering.}
  \label{fig:1d_case}
\end{figure}

Let us start by recalling a strategy to enhance \xfsai~\cite{Alsalti2023c} and later extend it towards \xamg.
Such a strategy introduces low-rank corrections to enable the use of \xspmm~but requires a discretisation consistent with the spatial symmetries.
Although boundary conditions need not be symmetric, they must be of the same type. Nevertheless, the final approach of \cref{ssec:MGR} will free this constraint, too.
Then, given an arbitrary mesh presenting a single reflection symmetry, we will order its grid points by first indexing the ones lying on one half and then those on the other.
Therefore, analogous to \cref{fig:1d_case}, if we impose the same local ordering (mirrored by the symmetry's hyperplane) to the resulting two subdomains, we ensure that all the scalar fields satisfy:
\begin{equation}\label{eq:fields_struct}
  x = 
    \begin{pmatrix}
      x_1 \\
      x_2
    \end{pmatrix} \in \R^n,
\end{equation}
where $n$ stands for the mesh size and $x_1,x_2\in\R^{n/2}$ for $x$'s restriction to each of the subdomains.
Mirrored grid points are in the same position within the subvectors, and discrete versions of virtually all partial differential operators satisfy the following block structure:
\begin{equation}\label{eq:genstruct_1sym}
  \H = 
    \begin{pmatrix}
      \H_{1,1} & \H_{1,2} \\
      \H_{2,1} & \H_{2,2}
    \end{pmatrix} \in \R^{n \times n},
\end{equation}
where $\H_{i,j}\in\R^{n/2 \times n/2}$ accounts for the couplings between the $i$th and $j$th subdomains.
As long as $\H$ only depends on geometric quantities (which is typically the case) or on material properties respecting the symmetries, given that both subdomains are identical and thanks to the mirrored ordering, we have that:
\begin{equation}\label{eq:prop_1sym}
  \H_{1,1}=\H_{2,2} \text{ and } \H_{1,2}=\H_{2,1},
\end{equation}
and, by denoting $\H_i \equiv \H_{1,i}$, we can rewrite \cref{eq:genstruct_1sym} as:
\begin{equation}\label{eq:struct_1sym}
  \H = 
    \begin{pmatrix}
      \H_1 & \H_2 \\
      \H_2 & \H_1
    \end{pmatrix}.
\end{equation}

The procedure above can be applied recursively to exploit an arbitrary number of reflection symmetries, $s$.
For instance, taking advantage of $s=2$ symmetries results in 4 mirrored subdomains and, analogously to \cref{eq:struct_1sym}, virtually all discrete operators satisfy the following:
\begin{equation}\label{eq:struct_2sym}
  \H = 
    \begin{pmatrix}
      \H_1 & \H_2 & \H_3 & \H_4 \\
      \H_2 & \H_1 & \H_4 & \H_3 \\
      \H_3 & \H_4 & \H_1 & \H_2 \\
      \H_4 & \H_3 & \H_2 & \H_1 
    \end{pmatrix},
\end{equation}
where $\H_i\in\R^{n/4 \times n/4}$ contains the couplings between the first and $i$th subdomains.

Thanks to the discretisation presented above, exploiting $s$ reflection symmetries allows meshing a $1/2^s$ fraction of the entire domain, henceforth named \emph{base mesh}.
\Cref{fig:dompart} depicts a schematic representation of the \emph{base mesh} on an arbitrary domain with two symmetries.
Then, instead of building the entire operators, $\H\in\R^{n \times n}$, it is only needed to build the \emph{base mesh}'s couplings with itself, $\H_1\in\R^{n/2^s \times n/2^s}$, and with its $2^s-1$ mirrored counterparts,  $\H_2, \dots, \H_{2^s}\in\R^{n/2^s \times n/2^s}$.
As a result, both the setup and memory footprint of the matrices are reduced by a factor of $2^s$~\cite{Alsalti2023b}.

Furthermore, while the sparsity pattern of $\H_1$ matches that of the actual operator built upon the \emph{base mesh}, the outer-subdomain couplings, $\H_2, \dots, \H_{2^s}$, have very few non-zero entries (if any), making the following splitting very advantageous:
\begin{equation}\label{eq:split_s}
  \H = \Id_{2^s} \otimes \Hinn + \Hout,
\end{equation}
where $\Hinn \coloneqq \H_1\in\R^{n/2^s \times n/2^s}$ and $\Hout \coloneqq \H - \Id_{2^s} \otimes \H_1\in\R^{n \times n}$.

The splitting of \cref{eq:split_s} allows the standard \xspmv~by $\H$ to be replaced with a specialised version of the more compute-intensive \xspmm.
The resulting kernel is a fusion of an \xspmm~by $\Hinn$, an \xspmv~by $\Hout$ and a linear combination of vectors, and has been employed throughout the simulations to make up to 5x faster their matrix multiplications~\cite{Alsalti2023b}.
In order to elucidate the computational advantages of \xspmm, let us note that the standard approach for applying $\Id_{2^s} \otimes \Hinn$ to a vector is through an \xspmv~call, which performs the following operation:
\begin{equation}\label{eq:spmv}
  y =
  \begin{pmatrix}
    \Hinn & & \\
     & \ddots & \\
      & & \Hinn
  \end{pmatrix}
  \begin{pmatrix}
    x_1 \\
    \vdots \\
    x_{2^s}
  \end{pmatrix} \in \Rn.
\end{equation}
Then, replacing \xspmv~with \xspmm~algebraically corresponds to:
\begin{equation}\label{eq:spmm}
  (y_1 \dots y_{2^s}) = \Hinn (x_1 \dots x_{2^s}) \in \R^{n/2^s \times 2^s}.
\end{equation}
The fact that \xspmm~reads $\Hinn$ $2^s$ fewer times makes its arithmetic intensity considerably higher and, since \xspmv~and \xspmm~are generally memory-bound kernels, this increase translates into significant speed-ups.

Apart from accelerating matrix multiplications, symmetries can be further harnessed to decompose Poisson's equation into a set of decoupled subsystems.
For the sake of clarity, let us recall the single-symmetry case of \cref{fig:1d_case}.
Following \cref{eq:split_s}, the discrete Laplacian can be written as:
\begin{equation}\label{eq:Lap_1sym}
  \A = 
    \begin{pmatrix}
      \A_1 & \A_2 \\
      \A_2 & \A_1
    \end{pmatrix} \in \R^{n \times n},
\end{equation}
and, by denoting the identity matrix with $\Id_k\in\R^{k \times k}$, we can define the following change-of-basis:
\begin{equation}\label{eq:P_s1}
  \P_1 \coloneqq \frac{1}{\sqrt{2}}
  \left(\begin{array}{c c}
    \Id_{\sfrac{n}{2}} & \Id_{\sfrac{n}{2}} \\
    \Id_{\sfrac{n}{2}} & -\Id_{\sfrac{n}{2}}
  \end{array}\right) \in \R^{n \times n},
\end{equation}
which transforms $\A$ into:
\begin{equation}\label{eq:hL_s1}
  \P_1 \A \P_1^{-1} =
  \left(\begin{array}{c c}
    \A_1 + \A_2 & \\
    & \A_1 - \A_2
  \end{array}\right).
\end{equation}

As shown in \cite{Alsalti2023a}, the block diagonalisation above can be generalised to an arbitrary number of reflection symmetries, $s$, by defining the following change-of-basis:
\begin{equation}\label{eq:P_s}
  \Ps \coloneqq \prod_{i=1}^s
  \left( \Id_{2^{i-1}} \otimes
  \frac{1}{\sqrt{2}}
  \begin{pmatrix}
    1 & 1 \\
    1 & -1
  \end{pmatrix} \otimes
  \Id_{\sfrac{n}{2^i}} \right) \in \Rnn,
\end{equation}
which satisfies $\Ps^{-1} = \Ps$ and transforms the discrete Laplacian into $2^s$ subsystems:
\begin{equation}\label{eq:hL_s}
  \hA \coloneqq \P_s \A \P_s^{-1} =
  \begin{pmatrix}
    \hA_1 &  &  \\
     & \ddots & \\
     & & \hA_{2^s}
  \end{pmatrix}.
\end{equation}
Then, similarly to \cref{eq:hL_s1}, $\hA$ can be split as follows:
\begin{equation}\label{eq:hL_splitting}
  \hA =
  \Id_{2^s} \otimes \Ainn +
  \begin{pmatrix}
    \A_{\text{out},1} &  &  \\
     & \ddots &  \\
     &  & \A_{\text{out},2^s}
  \end{pmatrix},
\end{equation}
making the Poisson solver of \cref{alg:symsolver} compatible with \xspmm.
Additionally, the decoupled solution of $\hA$'s subsystems in line \ref{alg:line:solve} makes Krylov subspace methods converge faster~\cite{Alsalti2023a}.

\begin{algorithm}
\caption{Poisson solver exploiting $s$ reflection symmetries}\label{alg:symsolver}
\begin{algorithmic}[1]
\Require $\hA_1,\dots,\hA_{2^s}$, $\Ps$ and $b \in \range(\A) \subseteq \Rn$
\Procedure{Solve}{$b$}\label{alg:line:proc_solve}
  \State {Transform forward: $\hat{b} = \Ps b$}\label{alg:line:ftrans}
  \State{Decoupled solution of $\hA_i\hat{x}_i = \hat{b}_i \;\; \forall i\in\{ 1,\dots, 2^s\}$}\label{alg:line:solve}
  \State {Transform backward: $x = \Ps \hat{x}$}\label{alg:line:btrans}
  \State \Return $x$
\EndProcedure
\end{algorithmic}
\end{algorithm}

\subsection{Low-rank corrections for factorable preconditioners}\label{ssec:lrc_sym}

In this section, we aim to develop low-rank corrections for enhancing factorable preconditioners, specifically focusing on making \xfsai~compatible with \xspmm.
Although initially introduced in \cite{Alsalti2023c}, our current objective is to extend this strategy towards non-factorable preconditioners like \xamg.

The idea of applying low-rank corrections arises from the close similarity between $\hA$'s subsystems.
Indeed, in \cref{eq:hL_splitting}, all the outer-couplings are substantially sparser than the inner and, as discussed in \cite{Alsalti2023c}, assuming a $d$-dimensional domain with a similar scale in all directions, we have that $\rank(\Aouti) = \mathcal{O}(n^{(d-1)/d})$, whereas $\rank(\Ainn) = \mathcal{O}(n)$.
Hence, for all $i \in \{1,\dots,2^s\}$:
\begin{equation}\label{eq:rank}
  \rank\left(\Aouti\right) \ll \rank\left(\Ainn\right),
\end{equation}
and it makes sense to introduce another level of approximation to \xfsai~by assuming that each of $\hA$'s subsystems satisfies:
\begin{equation}\label{eq:Ainn}
  \hA_i = \Ainn + \Aouti \simeq \Ainn.
\end{equation}

In the context of preconditioning linear systems, much work has recently been devoted to low-rank matrix representations~\cite{Li2013b,Ghysels2016,Li2017a,Franceschini2018,Xi2016}.
Let us recall the following result~\cite{Li2017a}.

\begin{theorem}\label{thm:correction}
  Given the two SPD matrices $A$ and $B$, let $L$ be the lower Cholesky factor of $B$, \ie~$B = LL^T$.
  Then, given $Y \coloneqq ( \Id - L^{-1} A L^{-T} )$, the following holds:
  \begin{equation*}
    A^{-1} = B^{-1} + L^{-T}V \Sigma \left( \Id - \Sigma \right)^{-1} V^TL^{-1},
  \end{equation*}
  where $Y=Y^T$ and $Y=V \Sigma V^T$ is the eigendecomposition of $Y$.
\end{theorem}

\begin{proof}
  From the definition of $Y$ it follows that $(\Id - Y)^{-1} = L^T A^{-1} L$ and, after some straightforward calculations:
  \begin{equation*}
    A^{-1} = B^{-1} + L^{-T} Y(\Id - Y)^{-1} L^{-1}.
  \end{equation*}
  Let $V \Sigma V^T$ be the eigendecomposition of $Y$. Then:
  \begin{equation*}
    (\Id - Y)^{-1} = V (\Id - \Sigma)^{-1} V^T
  \end{equation*}
  and, therefore, $A^{-1} = B^{-1} + L^{-T}V \Sigma \left( \Id - \Sigma \right)^{-1} V^TL^{-1}$.
\end{proof}

At this point, let us consider the \xfsai~of $\Ainn$.
It provides an approximation to the inverse of $\Ainn$'s lower Cholesky factor, $\Ginn \simeq L_{\text{inn}}^{-1}$, ensuring that:
\begin{equation}\label{eq:Ginn}
  \Ginn^T\Ginn \simeq \Ainn^{-1}.
\end{equation}

Then, for each subsystem $\hAi$, we can define the following auxiliary matrix:
\begin{equation}\label{eq:Ydef}
  Y \coloneqq \Id_{\sfrac{n}{2^s}} - \Ginn \hAi \Ginn^T \in \R^{n/2^s \times n/2^s},
\end{equation}
and, by virtue of \cref{thm:correction}, we have that:
\begin{equation}\label{eq:FRC}
  \hAi^{-1} = \Ginn^T\Ginn + \Ginn^T V \Sigma \left( \Id_{\sfrac{n}{2^s}} - \Sigma \right)^{-1} V^T \Ginn.
\end{equation}

The fact that $Y$'s eigendecomposition is dense prevents applying the ``full-rank'' correction of \cref{eq:FRC}.
However, thanks to \cref{eq:rank}, we can expect $Y$ to have a high data sparsity, \ie~its action is well represented by a low-rank approximation that only accounts for its $k$ most \emph{relevant} eigenpairs:
\begin{equation}\label{eq:Y_eigen}
  Y \simeq V_k \Sigma_k V_k^T,
\end{equation}
where $V_k \in \R^{n/2^s \times k}$ and $\Sigma_k \in \R^{k \times k}$ yield the following low-rank correction:
\begin{equation}\label{eq:LRC}
  \hAi^{-1} \simeq \Ginn^T\Ginn +  Z_k \Theta_k Z_k^T,
\end{equation}
with $Z_k \coloneqq \Ginn^TV_k \in \R^{n/2^s \times k}$ and $\Theta_k \coloneqq \Sigma_k (\Id_k - \Sigma_k)^{-1}\in \R^{k \times k}$.

As for selecting the most relevant eigenpairs, it is enough to remark that $Y$ measures how far each preconditioned subsystem, $\Ginn \hAi \Ginn^T$, is from the identity matrix.
Then, given the harmful effect that small eigenvalues have in the \xpcg~convergence~\cite{VanderSluis1986}, $Y$'s most effective eigenpairs are those associated with the smallest eigenvalues of $X \coloneqq \Ginn \hAi \Ginn^T$.
Hence, by computing a low-rank approximation of  $X$:
\begin{equation}\label{eq:X_eigen}
  X \simeq U_k \Lambda_k U_k^T,
\end{equation}
we can obtain the following truncated eigendecomposition of $Y$:
\begin{equation}\label{eq:YX_equiv}
  Y \simeq U_k (\Id_k - \Lambda_k) U_k^T.
\end{equation}
It can be shown experimentally that rough and cost-effective approximations of $U_k$ and $\Lambda_k$ suffice for preconditioning purposes.
Then, applying the above procedure to each of the $2^s$ subsystems leads to the following low-rank corrected \xfsai, henceforth denoted as \xlrcfsai:
\begin{equation}\label{eq:LRCFSAI}
  \Id_{2^s} \otimes \Ginn^T\Ginn + 
  \begin{pmatrix}
    Z_{k,1} \Theta_{k,1} Z_{k,1}^T & & \\
     & \ddots & \\
     & & Z_{k,2^s} \Theta_{k,2^s} Z_{k,2^s}^T
  \end{pmatrix}.
\end{equation}

\xlrcfsai~is compatible with \xspmm~and, thanks to reusing the same \xfsai~on all the subsystems, grants savings in the preconditioner's memory requirements and setup costs by a factor of $2^s$.
Of course, this comes at the price of using lower quality approximations, given that $\Ginn$ does not account for the outer-couplings.
However, introducing low-rank corrections proved very effective despite representing a relatively low overhead \cite{Alsalti2023c}.
It is remarkable that applying a rank-$k$ correction to each of the $2^s$ subsystems separately corresponds to applying a rank-$(2^sk)$ correction on the global system, $\hA$.

In order to illustrate the behaviour of \xlrcfsai~and, most especially, to assess the quality of the \xamg~variants that we will develop in the following sections, let us consider the following Poisson's equation with homogeneous Neumann boundary conditions:
\begin{equation}\label{eq:modelproblem}
	\begin{aligned}
	  -\frac{\partial^2 u}{\partial x^2} -\frac{\partial^2 u}{\partial y^2} -\frac{\partial^2 u}{\partial z^2} &= f \text{ in } \Omega \\
		\frac{\partial u}{\partial n} &= 0 \text{ on } \partial\Omega
	\end{aligned}
\end{equation}
where $f$ is a random field satisfying the compatibility condition $\iiint_\Omega f\;dV=0$.
The domain considered for the model problem is the unit cube, discretised using a standard 7-point stencil and the following hyperbolic stretching at the walls:
\begin{equation}\label{eq:tanh}
  x_i = \frac{1}{2}\left(1 + \frac{\tanh\left(\gamma_x \left(2\frac{(i-1)}{n_x}-1\right)\right)}{\tanh \left(\gamma_x\right)}\right) \;\; \forall i \in \{ 1, \dots, n_x + 1 \},
\end{equation}
analogously applied in the $y$- and $z$-directions using $\gamma_x=\gamma_y=\gamma_z=1.5$.

\Cref{tab:lrcfsai} shows the convergence of \xlrcfsai{} on the model problem of \cref{eq:modelproblem}.
The results correspond to a MATLAB implementation using an FSAI with up to 50 nonzeros per row and a pre-filtering tolerance $\tau=0.01$~\cite{JanFerSarGam15}.
Clearly, low-rank corrections outweigh the fact of ignoring the outer-couplings, actually making \xlrcfsai{} converge faster than the standard \xfsai.
As a result, we proceeded with its parallel implementation, which is out of the scope of this work but was thoroughly reviewed in \cite{Alsalti2023c}.

\begin{table}[htbp]
{\footnotesize
  \caption{PCG + \protect\xlrcfsai~results on the model problem with $n=64^3$.}\label{tab:lrcfsai}
  \begin{center}
  	\begin{tabularx}{0.7\textwidth}{ |c| *{4}{Y|} }
	  \cline{2-5}
  	   \multicolumn{1}{c|}{} 
	   & \multicolumn{4}{c|}{iterations}\\
  	\hline
	   preconditioner & $s=0$ & $s=1$ & $s=2$ & $s=3$ \\
  	\hline
     FSAI        & 291 & 220 & 173 & 121 \\
     LRCFSAI(0)  & 291 & 241 & 193 & 149 \\
     LRCFSAI(1)  & 228 & 194 & 146 & 106 \\
     LRCFSAI(2)  & 228 & 172 & 124 & 92  \\
     LRCFSAI(4)  & 186 & 136 & 107 & 83  \\
     LRCFSAI(8)  & 131 & 111 & 81  & 61  \\
     LRCFSAI(16) & 113 & 82  & 63  & 50  \\
  	\hline
	  \end{tabularx}
  \end{center}
}
\end{table}

\subsection{Low-rank corrections for non-factorable preconditioners}\label{ssec:lrc_nonsym}

At this point, we can extend the strategy presented in \cref{ssec:lrc_sym} towards \xamg, which is probably the most powerful preconditioner for Poisson's equation.
The effectiveness of \xamg~is given by the complementary roles played by the smoother, which is responsible for damping high-frequency error components, and the coarse-grid correction, which in turn reduces low-frequency modes.
Large problems require a progressive coarsening into a hierarchy of smaller and smaller grids.
These grids are created using the concept of strength of connection, which measures the likelihood that the smooth error components on two adjacent nodes have similar values.
There are a number of beautiful books on the subject, \eg~\cite{Tro01,Vas08}.

Bearing in mind that the application of \xamg~relies on matrix multiplications, we aim to accelerate it by replacing \xspmv~with \xspmm.
However, \xamg~is not explicitly factorable, and we need to develop alternative low-rank corrections.
With this aim, let us recall the following result~\cite{Li2017a},
in which we replace SVD with eigendecomposition, as it facilitated the selection of the correction vectors.

\begin{theorem}\label{thm:correction_nonfact}
  Given the two SPD matrices $A$ and $B$, let us define the auxiliary matrix $Y \coloneqq ( \Id - B^{-1} A )$.
  Then, the following holds:
  \begin{equation*}
    A^{-1} = \left( \Id + U \left( \Sigma^{-1} - V^TU \right)^{-1} V^T \right) B^{-1},
  \end{equation*}
  where $Y\neq Y^T$ and $Y=U \Sigma V^T$ is the eigendecomposition of $Y$.
\end{theorem}

\begin{proof}
  From the definition of $Y$, it follows that $A^{-1} = (\Id - Y)^{-1} B^{-1}$.
  Then, applying the Sherman-Morrison-Woodbury formula, we have that:
  \begin{equation*}
    (\Id - Y)^{-1} = (\Id - U \Sigma V^T)^{-1} = \Id + U (\Sigma^{-1} - V^TU)^{-1} V^T,
  \end{equation*}
  and combining both equations completes the proof.
\end{proof}

As we did to derive \xlrcfsai, let us start by considering the following \xamg{} approximation of $\Ainn^{-1}$:
\begin{equation}\label{eq:AMGinn}
  \AMGinn \simeq \Ainn^{-1}.
\end{equation}
Then, for each subsystem $\hAi$, we can define the following auxiliary matrix:
\begin{equation}\label{eq:Ydef_AMG}
  Y \coloneqq \Id_{\sfrac{n}{2^s}} - \AMGinn \hAi \in \R^{n/2^s \times n/2^s},
\end{equation}
and, by virtue of \cref{thm:correction_nonfact}, we have that:
\begin{equation}\label{eq:FRC_AMG}
  \hAi^{-1} = \left( \Id_{\sfrac{n}{2^s}} + U \left( \Sigma^{-1} - V^TU \right)^{-1} V^T \right) \AMGinn.
\end{equation}
Once again, thanks to \cref{eq:rank}, we can expect $Y$ to be well represented by a truncated eigendecomposition only accounting for its $k$ most \textit{relevant} eigenvectors:
\begin{equation}\label{eq:Y_eigen_AMG}
  Y \simeq U_k \Sigma_k V_k^T,
\end{equation}
where $U_k, V_k \in \R^{n/2^s \times k}$ and $\Sigma_k \in \R^{k \times k}$.
Differently to \cref{eq:Y_eigen} in \cref{ssec:lrc_sym}, $Y$'s nonsymmetry makes it require both the right- and left-eigenvectors, yielding the following low-rank correction:
\begin{equation}\label{eq:LRC_AMG}
  \hAi^{-1} \simeq
	  \left( \Id_{\sfrac{n}{2^s}} + U_k ( \Sigma_k^{-1} - V_k^TU_k )^{-1} V_k^T \right) \AMGinn =
  	\left( \Id_{\sfrac{n}{2^s}} + U_k \Theta_k V_k^T \right) \AMGinn,
\end{equation}
where we considered biorthonormal bases for the right- and left-eigenvectors, \ie{} $V_k^TU_k = \Id_k$, and defined $\Theta_k \coloneqq ( \Sigma_k^{-1} - \Id_k )^{-1} \in \R^{k \times k}$.

Then, applying the above procedure to each of the $2^s$ subsystems results in the following low-rank corrected \xamg, henceforth denoted as \xlrcamg:
\begin{equation}\label{eq:LRCAMG}
  \left[ \Id_n + 
  \begin{pmatrix}
    U_{k,1} \Theta_{k,1} V_{k,1}^T & &  \\
     & \ddots & \\
      & & U_{k,2^s} \Theta_{k,2^s} V_{k,2^s}^T
  \end{pmatrix} \right]
  \left(\Id_{2^s} \otimes \AMGinn\right).
\end{equation}

When it comes to selecting the most relevant eigenvectors, let us note that, as in \cref{ssec:lrc_sym}, $Y$ measures how far it is each preconditioned subsystem, $\AMGinn \hAi$, from the identity matrix.
Then, given the harmful effect of small eigenvalues in the convergence of Krylov subspace methods, $Y$'s most effective eigenvectors are those associated with the smallest eigenvalues of $X \coloneqq \AMGinn \hAi$.
Additionally, by computing the eigendecomposition of $X$:
\begin{equation}\label{eq:X_eigen_AMG}
  X = U \Lambda V^T,
\end{equation}
we can obtain the following eigendecomposition of $Y$:
\begin{equation}\label{eq:YX_equiv_AMG}
  Y = U (\Id - \Lambda) V^T,
\end{equation}
which gives $\Sigma = \Id - \Lambda$.
Hence, the smallest, and harmful, eigenvalues of $X$ correspond to the largest positive eigenvalues of $Y$, whereas large (in magnitude) negative eigenvalues of $Y$ correspond to the largest, and harmless, eigenvalues of $X$.

\xlrcamg~is compatible with \xspmm~in the application of $\AMGinn$.
However, this comes at the price of introducing another level of approximation by ignoring $\hA$'s outer-couplings.
While low-rank corrections proved very effective on \xlrcfsai, even accelerating its convergence, (relatively) low-rank perturbations had a critical impact on \xamg.
Indeed, \cref{tab:lrcamg} summarises the results obtained with \xlrcamg~on the model problem of \cref{eq:modelproblem}.
The results correspond to a MATLAB implementation using an AMG with a PMIS coarsening~\cite{Sterck2006}, Extended+I interpolation~\cite{DeSFalNolYan08} and an FSAI smoother with 15 nonzeros per row.
Unfortunately, low-rank corrections cannot restore the effectiveness of \xamg, which excels in removing error components on the lower part of the spectrum with high accuracy, something that low-rank correction cannot match.
Due to its unsatisfactory performance, we abandoned its parallel implementation and explored the better alternatives of \cref{sec:innifc}.

\begin{table}[htbp]
{\footnotesize
  \caption{GMRES + \protect\xlrcamg~results on the model problem with $n=64^3$.}\label{tab:lrcamg}
  \begin{center}
  	\begin{tabularx}{0.7\textwidth}{ |c| *{4}{Y|} }
	  \cline{2-5}
  	   \multicolumn{1}{c|}{} 
	   & \multicolumn{4}{c|}{iterations}\\
  	\hline
	   preconditioner & $s=0$ & $s=1$ & $s=2$ & $s=3$ \\
  	\hline
      AMG        & 6 & 6  & 6  & 6  \\
      LRCAMG(0)  & 6 & 23 & 29 & 29 \\
      LRCAMG(1)  & 9 & 23 & 26 & 24 \\
      LRCAMG(2)  & 9 & 23 & 22 & 24 \\
      LRCAMG(4)  & 9 & 19 & 21 & 20 \\
      LRCAMG(8)  & 9 & 17 & 18 & 18 \\
      LRCAMG(16) & 9 & 15 & 16 & 15 \\
  	\hline
	  \end{tabularx}
  \end{center}
}
\end{table}

\section{Inner-interface spatial discretisation}\label{sec:innifc}

\begin{figure}[tbp]
	\centering
  \begin{subfigure}{\linewidth}
	  \centering
  	\includegraphics[trim = 5.5cm 9.6cm 18.9cm 12.5cm,clip,width=\linewidth]{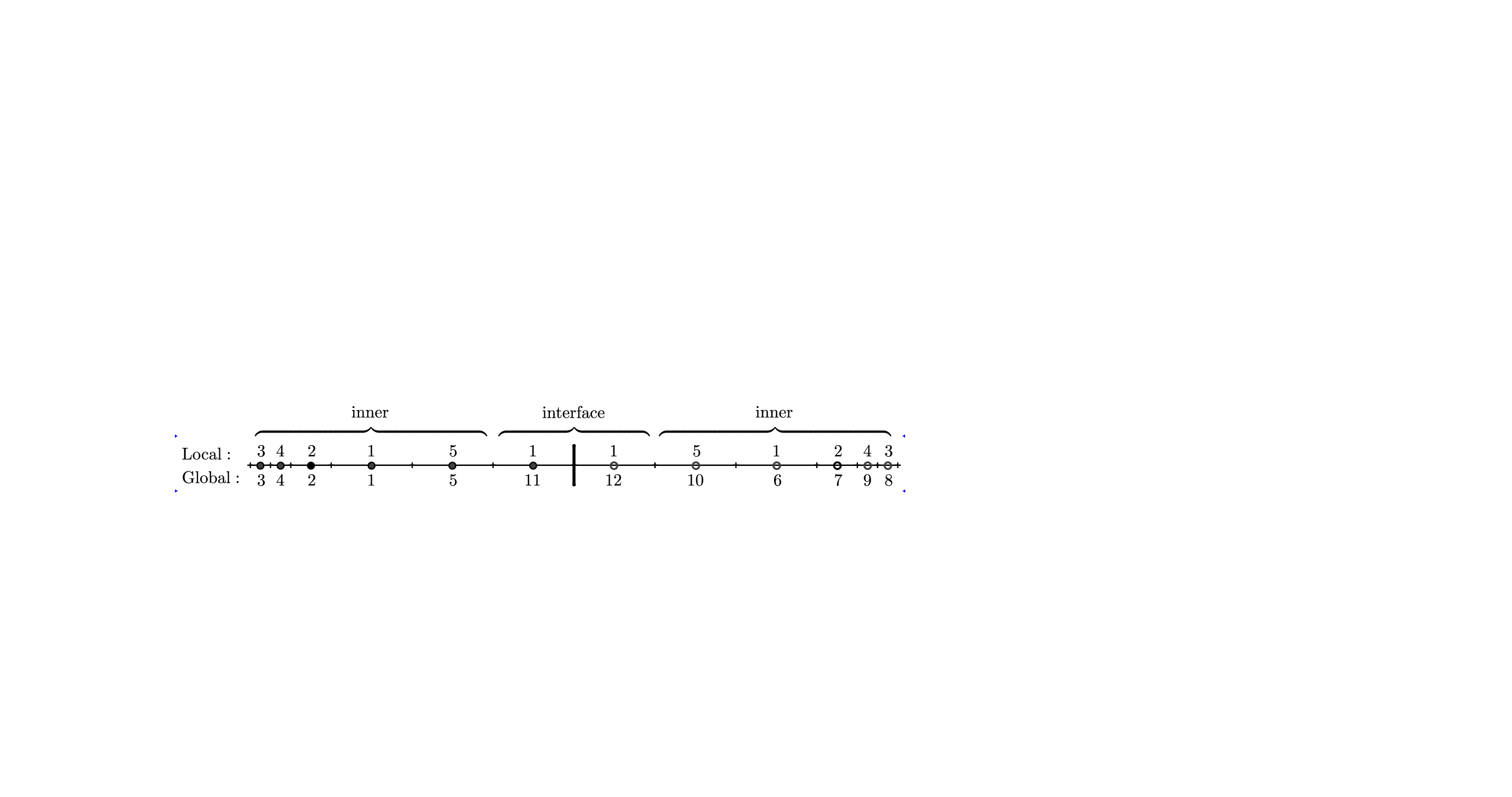}
	  \caption{Mirrored geometry}
    \label{fig:1d_case_innifc_sym}
  \end{subfigure}
  \begin{subfigure}{\linewidth}
	  \centering
  	\includegraphics[trim = 5.5cm 9.6cm 18.9cm 12.5cm,clip,width=\linewidth]{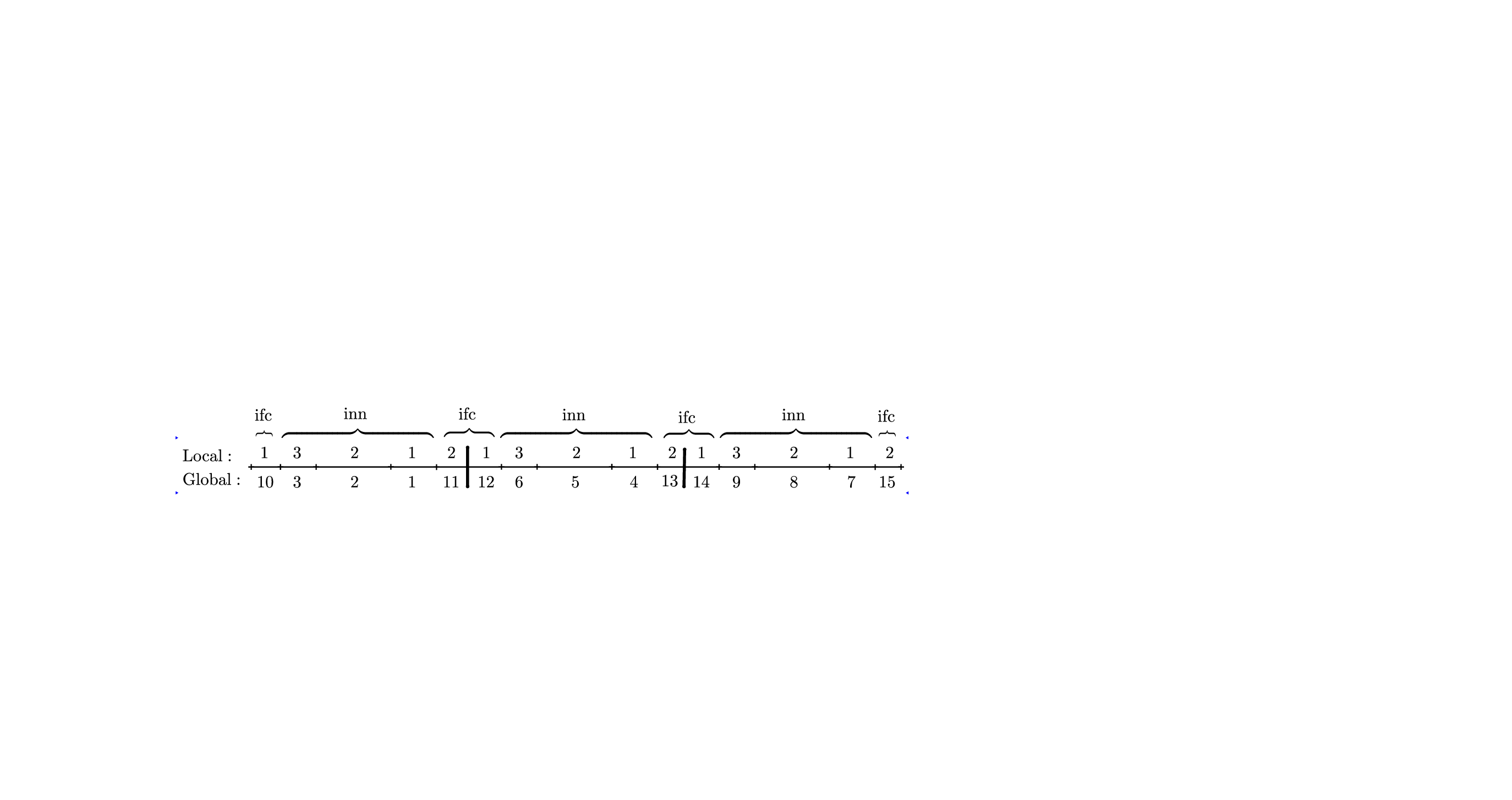}
    \caption{Repeated geometry}
    \label{fig:1d_case_innifc_rep}
  \end{subfigure}
	\caption{1D meshes with a random inner-interface ordering.}
	\label{fig:1d_case_innifc}
\end{figure}

This section explores two new approaches to accelerate \xamg~by making it compatible with \xspmm.
Unlike \xlrcfsai{} and \xlrcamg, the preconditioners developed here rely on a different ordering of unknowns, henceforth denoted as inner-interface ordering.
To define it, let us recall the single symmetry case of \cref{fig:1d_case_innifc_sym} and classify all the grid points into interface: those coupled with other subdomain unknowns, and inner: those that are not.
Then, the inner-interface ordering arises from applying the mirrored ordering of \cref{fig:1d_case} to the inner unknowns first and the interface ones afterwards.
Remarkably enough, AMGR, the preconditioner derived in \cref{ssec:MGR}, relaxes this constraint, allowing any ordering of the interface unknowns and, most notably, freeing all constraints on the boundary conditions.
Hence, AMGR, our final proposal, will not only succeed in outperforming the standard \xamg~but also apply to a broader range of industrial applications, including, for instance, repeated geometries like the illustrative \cref{fig:1d_case_innifc_rep}.

As a result of using the inner-interface ordering, the discrete Laplacian reads:
\begin{equation}\label{eq:Linnifc}
  \A = 
    \begin{pmatrix}
      \bK & \bB \\
      \bBt & \bC
    \end{pmatrix} \in \R^{n \times n},
\end{equation}
where $\bK \in\R^{\ninn \times \ninn}$, $\bB \in\R^{\ninn \times \nifc}$, and $\bC \in\R^{\nifc \times \nifc}$.
Additionally, thanks to the mirrored ordering imposed on the inner and interface grid points, we have that $\bK = \Id_{2^s} \otimes \K$ and $\bB = \Id_{2^s} \otimes \B$.
On the other hand, $\bC$ satisfies \cref{eq:split_s}, and consequently, \xspmm~can be exploited to accelerate all the submatrix products while reducing their setup costs and memory footprint.

\subsection{Schur complement-based AMG}\label{ssec:SchurAMG}

Given the poor results obtained with \xlrcamg, our aim now is to exploit the structure of $\bK$, $\bB$ and $\bC$ to accelerate the \xamg{} approximation of $\A$, $\AMGA\simeq\A^{-1}$.
By virtue of \cref{eq:Linnifc}, the inverse of $\A$'s LDU factorisation reads:
\begin{equation}\label{eq:L_domdec}
  \A^{-1} = 
  \begin{pmatrix}
    \Id & -\bK^{-1}\bB  \\
        & \Id
  \end{pmatrix}
  \begin{pmatrix}
    \bK^{-1} &          \\
             & \bS^{-1}
  \end{pmatrix}
  \begin{pmatrix}
    \Id           &     \\
    -\bBt\bK^{-1} & \Id
  \end{pmatrix},
\end{equation}
where $\bS \coloneqq \bC - \bBt \bK^{-1} \bB$ is the Schur complement of $\bK$.
Then, we can derive valid preconditioners by seeking approximations of $\bK^{-1}$ and $\bS^{-1}$.
In our case, given that $\bK^{-1} = \Id_{2^s} \otimes \K^{-1}$, we will approximate $\bK^{-1}$ through the AMG of $\K$, $\AMGK \simeq \K^{-1}$.
This has two major advantages: on the one hand, the setup costs and memory footprint of $\AMGK$ are $2^s$ times smaller than those of $\AMGbK$.
On the other, its application is compatible with \xspmm.

Regarding the approximation of $\bS^{-1}$, the better it is, the closer the action of our Schur-based \xamg~will be to $\AMGA$.
We have that:
\begin{equation}\label{eq:bS}
  \bS \coloneqq \bC - \bBt \bK^{-1} \bB
	    = \bC - \Id_{2^s} \otimes \Bt \K^{-1} \B,
\end{equation}
and $\bS$ containing $\K^{-1}$ prevents us from using \xamg~to approximate $\bS^{-1}$.
However, we can obtain a straightforward approximation of $\bS^{-1}$ using \xfsai.
With this aim, let $L_{\K}$ be the lower Cholesky factor of $\K$, and $\GK \simeq L_{\K}^{-1}$ its lower \xfsai~factor.
Then, we have that $\GK^T \GK \simeq \K^{-1}$, and we can define:
\begin{equation}\label{eq:bSG}
  \bC - \Id_{2^s} \otimes \Bt \GK^T \GK \B \simeq \bS,
\end{equation}
whose inverse can be approximated through another \xfsai, leading to an approximation that we will denote as $\bS_{\G}^{-1} \simeq \bS^{-1}$.

As the model problem will demonstrate, such an approximation of $\bS^{-1}$ is too rough, resulting in slow convergences.
Recalling the favourable results of subsection \ref{ssec:lrc_sym}, we will build an approximation of $\bS^{-1}$ based on \xlrcfsai.
For simplicity, let us assume a single reflection symmetry.
As before, applying a mirrored ordering to the interface grid points ensures that $\bC$ satisfies the following structure:
\begin{equation}\label{C_p1}
  \bC =
  \begin{pmatrix}
    \C_1 & \C_2 \\
    \C_2 & \C_1
  \end{pmatrix},
\end{equation}
where $\C_1,\C_2\in\R^{\nifc/2 \times \nifc/2}$ account for the interface-interface couplings of the \emph{base mesh} with itself and its mirroring, respectively.
Combining \cref{C_p1,eq:bS}, we have that:
\begin{equation}\label{S_p1}
  \bS = \begin{pmatrix}
		  	  \S_1 & \C_2 \\
			    \C_2 & \S_1
			  \end{pmatrix},
\end{equation}
where we defined $\S_1 \coloneqq \C_1 - \B^T \K^{-1} \B$.
Then, in order to compute $\bS^{-1}$  we can proceed as in \cref{eq:Lap_1sym,eq:P_s1,eq:hL_s1} to block diagonalise the Schur complement of $\bK$:
\begin{equation}\label{PSP_p1}
  \hS \coloneqq \P_1\bS\P_1^{-1} =
  \begin{pmatrix}
    \S_1 + \C_2 & \\
    & \S_1-\C_2
  \end{pmatrix}.
\end{equation}

After the excellent results obtained with \xlrcfsai{} and its decoupled corrections, it makes sense to mimic that strategy.
First, we will build a rough approximation of $\hS^{-1}$, similar to $\bS_{\G}^{-1}$ but based on an \xfsai~of $\S_1$ and, therefore, compatible with \xspmm.
Then, we will correct it through low-rank corrections to make it approach the (block diagonal version of the) most accurate Schur complement at hand:
\begin{equation}\label{eq:S_AMGA}
  \bS_{\text{AMG}} \coloneqq \bC - \Id_{2^s} \otimes \B^T \AMGK \B,
\end{equation}
but whose inverse is unavailable.

Similarly to \cref{eq:bSG}, it is meaningful to approximate $\S_1^{-1}$ through an \xfsai~of $\C_1 - \B^T \GK^T \GK \B \simeq \S_1$:
\begin{equation}\label{eq:invS_GA}
   \G_{\S_1}^T \G_{\S_1} \simeq \S_1^{-1}.
\end{equation}
Then, by recalling \cref{PSP_p1} and the fact that $\S_1$ is substantially denser than $\C_2$, we can build the sought approximation of $\hS^{-1}$ as follows:
\begin{equation}\label{eq:invS_GA_inn}
  \Id_2 \otimes \G_{\S_1}^T \G_{\S_1} \simeq \hS^{-1},
\end{equation}
which not only is compatible with \xspmm~but also allows invoking \cref{thm:correction} on each decoupled block:
\begin{equation}\label{eq:S_LRC_inn}
  \hS_\text{LRC}^{-1} \coloneqq 
  \Id_2 \otimes \G_{\S_1}^T \G_{\S_1} + 
  \begin{pmatrix}
    Z_k^\text{(1)} \Theta_k^\text{(1)} {Z_k^\text{(1)}}^T & \\
    & Z_k^\text{(2)} \Theta_k^\text{(2)} {Z_k^\text{(2)}}^T
  \end{pmatrix}.
\end{equation}

Remarkably enough, to introduce low-rank corrections that improve \xfsai, it is not necessary to leverage symmetries and block diagonalise $\bS$.
Instead, we could have invoked \cref{thm:correction} on $\bS_{\G}^{-1}$, correcting it to approach $\bS_{\text{AMG}}$ and giving rise to $\bS_\text{LRC}^{-1}\simeq \bS^{-1}$.
However, as discussed in \cref{ssec:lrc_sym} and demonstrated with the model problem, the larger the number of subdomains, the more effective the decoupled corrections.

Regardless of approximating $\bS^{-1}$ with $\bS_{\G}^{-1}$, $\bS_\text{LRC}^{-1}$, or $\Ps^{-1} \hS_\text{LRC}^{-1} \Ps$, \cref{eq:L_domdec} requires inverting $\bK$ twice, whereas $\AMGA$ is only applied once.
In order to avoid such an expensive overhead, we can define our Schur-based \xamg, henceforth denoted AMGS$(\bS^{-1}, \bK^{-1})$, as follows:
\begin{equation}\label{eq:M_domdec}
  \AMGS(\bS^{-1}, \bK^{-1}) \coloneqq \begin{pmatrix}
    \Id & -\bK^{-1}\bB  \\
        & \Id
  \end{pmatrix}
  \begin{pmatrix}
    \AMGbK & \\
           & \bS^{-1}
  \end{pmatrix}
  \begin{pmatrix}
    \Id          &     \\
    -\bBt \AMGbK & \Id
  \end{pmatrix},
\end{equation}
where we compacted the notation of $\AMGbK \coloneqq \Id_{2^s} \otimes \AMGK$ and replaced the extra application of $\AMGK$ with $\bK^{-1}$, which, for instance, could be the \xfsai~of $\bK$, $\G_{\bK}^T \G_{\bK} \coloneqq \Id_{2^s} \otimes \GK^T \GK$.
The main advantages of AMGS$(\bS^{-1}, \bK^{-1})$ are being its matrix multiplications compatible with \xspmm, together with the smaller memory footprint and setup costs that this allows.

\Cref{tab:shuramg} summarises the results obtained with AMGS on the model problem of \cref{eq:modelproblem}.
The results correspond to a MATLAB implementation using an AMG with a PMIS coarsening~\cite{Sterck2006}, Extended+I interpolation~\cite{DeSFalNolYan08} and an FSAI smoother with 15 nonzeros per row.
Conversely, the FSAI of $\K$ and of the Schur complement have a pre-filtering tolerance $\tau=0.01$~\cite{JanFerSarGam15}, and 100 and 50 nonzeros per row, respectively.
For completeness, three approximations of $\bS^{-1}$ are included in \cref{tab:shuramg}.
Namely, $\bS_{\G}^{-1}$, $\bS_\text{LRC}^{-1}$, and $\Ps^{-1} \hS_\text{LRC}^{-1} \Ps$.
The fact that AMGS$(\bS_{\G}^{-1},\AMGbK)$ converges almost twice slower than the other AMGS variants confirms the need for low-rank corrections.
Especially as $s$ grows, since $\nifc$ grows accordingly and convergence becomes more sensitive to the approximation of $\bS^{-1}$.
Due to this, the convergence of AMGS$(\bS_\text{LRC}^{-1}, \AMGbK)$ also deteriorates with $s$.
However, thanks to becoming more effective as the number of subdomains grows, decoupled corrections do not suffer from this effect, requiring AMGS$(\hS_\text{LRC}^{-1}, \AMGbK)$ roughly the same iterations regardless of $s$.
Note that AMGS$(\hS_\text{LRC}^{-1}, \AMGbK)$ does actually refer to AMGS$(\Ps^{-1} \hS_\text{LRC}^{-1} \Ps, \AMGbK)$.

All three variants above require applying $\AMGbK$ twice, which makes the computational costs of the resulting preconditioner exceedingly large.
However, as hinted earlier, we can replace the extra application of $\AMGbK$ with \xfsai{} in \cref{eq:M_domdec}.
\Cref{tab:shuramg} confirms that we can do this safely, and although it is much lighter, AMGS$(\hS_\text{LRC}^{-1}, \GbK)$ converges comparably to AMGS$(\hS_\text{LRC}^{-1}, \AMGbK)$.

Unfortunately, the combination of \xfsai~and low-rank corrections was not enough to substantially improve the performance of AMGS with respect to \xlrcamg{} in \cref{tab:shuramg}.
Indeed, it converged more than three times slower than \xamg, rendering it impossible to reduce the solution time regardless of the more compute-intensive matrix multiplications.
Consequently, we discarded its parallel implementation and developed the final alternative of \cref{ssec:MGR}, which converges comparably to AMG and offers greater flexibility than all previous alternatives.

\begin{table}[htbp]
{\footnotesize
  \caption{GMRES + AMGS results on the model problem with $n=64^3$.}\label{tab:shuramg}
  \begin{center}
  	\begin{tabularx}{0.7\textwidth}{ |c| *{4}{Y|} }
	  \cline{2-5}
  	   \multicolumn{1}{c|}{} 
	   & \multicolumn{4}{c|}{iterations}\\
  	\hline
	   preconditioner & $s=0$ & $s=1$ & $s=2$ & $s=3$ \\
  	\hline
    AMG                                & 6  & 6  & 6  & 6 \\
    AMGS$(\bS_{\G}^{-1},       \AMGbK)$ & 6 & 35 & 38 & 44 \\
    AMGS$(\bS_\text{LRC}^{-1}, \AMGbK)$ & 6 & 22 & 26 & 29 \\
    AMGS$(\hS_\text{LRC}^{-1}, \AMGbK)$ & 6 & 20 & 23 & 21 \\
    AMGS$(\hS_\text{LRC}^{-1}, \GbK)$   & 6 & 21 & 23 & 22 \\ 
  	\hline
	  \end{tabularx}
  \end{center}
}
\end{table}

\subsection{Multigrid reduction}\label{ssec:MGR}

As the model problem of \cref{eq:modelproblem} confirmed, the standard \xamg~algorithm cannot take advantage of spatial symmetries and needs to be appropriately adapted.
In particular, we will develop an \xamg~reduction framework, hereafter denoted as AMGR, that converges comparably to \xamg, is compatible with \xspmm, and leverages reflection, translational and rotational symmetries.
Remarkably enough, AMGR only uses \xspmm~on the inner unknowns, therefore allowing any interface ordering and, more importantly, freeing all constraints on the boundary conditions, which makes it apply to a significantly broader range of industrial applications without compromising its computational advantages, as $\ninn\gg\nifc$.

Ideally, to have the fastest possible coarsening, we should only identify as coarse the interface unknowns, \ie~the $\nifc$ unknowns associated with the interface block, $\bC$.
However, realistically sized problems make it impossible to accurately interpolate inner unknowns solely using interface values.
This is due to the large connection distance that may occur.
Indeed, assuming a $d$-dimensional domain with a similar scale in all directions and $\ninn$ inner unknowns, the maximum inner-interface distance will be of the order of $\sqrt[d]{\ninn}$.
Then, a 3D problem with about a million inner unknowns would result in distances of about 100 units, which exceeds the applicability of long-distance interpolation formulas such as Extended+I (ExtI)~\cite{DeSFalNolYan08} or Dynamic Pattern Least Squares (DPLS)~\cite{PalFraJan19}.

Consequently, to allow for an accurate interpolation, we need to convert some inner nodes into coarse.
In principle, we can apply any standard coarsening strategy to the inner part of the \textit{base mesh}.
That is, to choose a strength of connection measure, \eg~a classical strength measure, filtering the resulting adjacency graph to obtain $T$, and running a maximum independent set (MIS) algorithm on $T$ to select which inner nodes become coarse.
Such a procedure populates the interface with some of the \textit{base mesh} inner unknowns.
Then, to preserve the structure of $\bK = \Id_{n_b} \otimes \K$, we need their analogues in the rest of the subdomains to be set as coarse, too.
The problem is that, even if many connections are retained in $T$, the number of coarse nodes generally becomes too large, making the resulting operator's complexity impractical.

\Cref{alg:popul} summarises the aggressive coarsening strategy designed ad hoc for tackling this problem.
Instead of selecting independent nodes from $T$, we select independent nodes from $T^k$, for a small power $k$.
Using $T^k$ is equivalent to constructing an independent set by considering dependent on each other all nodes at a distance smaller or equal to $k$.
For simplicity, from now on, $\nifc$ and $\ninn$ will refer to the populated interface and the remaining inner unknowns.
That is, to the number of coarse and fine nodes arising from \cref{alg:popul}, respectively.
The same applies to the partitioning of \cref{eq:Linnifc} and $\bK$, $\bB$ and $\bC$.

\begin{algorithm}[htbp]
\caption{Populated coarsening at a maximum interpolation distance $k$}\label{alg:popul}
\begin{algorithmic}[1]
\Procedure{AMGR\_Coarsening}{$A$, $k$}
  \State Measure the strength-of-connection of $A$
  \State Filter the resulting adjacency graph, $T$
  \State Compute the symbolic power of $T$, $T^k$
  \State Compute a maximum independent set on $T^k$, MIS($T^k$)
  \State Set as coarse the unknowns $\text{MIS}(T^k) \;\dotcup\; \{\nifc+1,\dots,n\}$
\EndProcedure
\end{algorithmic}
\end{algorithm}

The coarsening of \cref{alg:popul} yields a prolongation with the following structure:
\begin{equation}\label{eq:amgreduction}
P \coloneqq
\begin{pmatrix}
\bar{W} \\
\Id_{\nifc}
\end{pmatrix} \in \R^{n \times \nifc},
\end{equation}
where $\bar{W} \in\R^{\ninn \times \nifc}$ and $\Id_{\nifc} \in\R^{\nifc \times \nifc}$.
Relatively large distances between coarse nodes are not an issue if long-distance interpolations are used.
ExtI interpolation is most effective when $k$ is limited to 2.
Conversely, DPLS allows for larger values of $k$ and coarser operators.
However, there is a trade-off between increasing the value of $k$ to have a more aggressive coarsening and keeping the interpolations accurate.
\Cref{sec:experiments} investigates the optimal configuration for AMGR, concluding that best results are obtained with $k=2$, DPLS and energy minimisation~\cite{JanFraSchOls23}.

In order to make AMGR exploit the computational advantages of \xspmm, the top-level smoother is defined as the following block FSAI:
\begin{equation}\label{eq:topsmoother}
  \begin{pmatrix}
    \GbK^T\GbK & \\
    & \GbC^T\GbC
  \end{pmatrix} \simeq \A^{-1},
\end{equation}
where $\GbK = \Id_{n_b} \otimes \GK$, $\GK^T\GK \simeq \K^{-1}$ and $\GbC^T\GbC \simeq \bC^{-1}$.
Numerical experiments confirmed that, for relaxation purposes, one can safely ignore $A$'s off-diagonal blocks in \cref{eq:topsmoother}, as it causes no significant convergence degradation.
Nevertheless, it allows the replacement of \xspmv{} with \xspmm{} and $\GbK$ with $\GK$, therefore reducing the memory footprint, setup and application costs of the top-level smoother.

Finally, the reduced operator resulting from the multigrid reduction reads:
\begin{equation}\label{eq:structAc}
  A_c \coloneqq
  P^T A
  P =
  \bW^T\bK\bW + \bW^T\bB + \bBt\bW + \bC \in \R^{\nifc \times \nifc},
\end{equation}
which is then inverted through a standard \xamg{}.
\Cref{alg:apply} summarises the application of AMGR, which aims to enhance the \xamg{} applied in line \ref{line:amgapply} by reducing the costs of its finest (and most computationally expensive) level.
Hence, the multigrid reduction only defines the first level of AMGR, and subsequent levels are purely algebraic.

\begin{algorithm}[htbp]
\caption{AMGR application in a V$(\nu_1,\nu_2)$-cycle}\label{alg:apply}
\begin{algorithmic}[1]
\Procedure{AMGR\_Apply}{$A$, $y$, $z$}
  \State Smooth $\nu_1$ times $A s = y$ starting from $0_n$
  \State Compute the residual $r = y - A s$
  \State Restrict the residual to the populated interface grid $r_c = P^T r$
  {\State \text{AMG\_Apply}$\left(A_c,r_c,h_c \right)$}\label{line:amgapply}
  \State Prolongate the correction to the entire grid $h = P h_c$
  \State Update $s = s + h$
  \State Smooth $\nu_2$ times $A z = y$ starting from $s$
\EndProcedure
\end{algorithmic}
\end{algorithm}

\Cref{tab:amgr} summarises the results obtained with AMGR on the model problem of \cref{eq:modelproblem}.
The results correspond to a MATLAB implementation of a multigrid reduction using $k=2$, DPLS and energy minimisation.
Conversely, the AMG applied to the reduced operator used a PMIS coarsening~\cite{Sterck2006}, Extended+I interpolation~\cite{DeSFalNolYan08} and an FSAI smoother with 15 nonzeros per row.
As it can be seen, adopting the aggressive coarsening that spatial symmetries induce does not harm convergence and allows for the computational advantages of replacing \xspmv{} with \xspmm{}.
Hence, given such promising results, we tackled its parallel implementation, which is discussed in \cref{sec:impl}, and evaluated its performance on industrial \xcfd~applications in \cref{sec:experiments}.

\begin{table}[htbp]
{\footnotesize
  \caption{PCG + AMGR results on the model problem with $n=64^3$.}\label{tab:amgr}
  \begin{center}
  	\begin{tabularx}{0.7\textwidth}{ |c| *{4}{Y|} }
	  \cline{2-5}
  	   \multicolumn{1}{c|}{} 
	   & \multicolumn{4}{c|}{iterations}\\
  	\hline
	   preconditioner & $s=0$ & $s=1$ & $s=2$ & $s=3$ \\
  	\hline
	   AMG  & 6 & 6 & 6 & 6 \\
	   AMGR & 6 & 9 & 9 & 9 \\
  	\hline
	  \end{tabularx}
  \end{center}
}
\end{table}

\section{Practical implementation}\label{sec:impl}

Discretising complex geometries becomes simpler thanks to exploiting spatial symmetries.
Indeed, the strategies presented only require meshing the \emph{base mesh}, which, assuming $n_b$ subdomains, corresponds to a $1/n_b$ fraction of the entire domain.
Then, the implementation expands the \emph{base mesh} by imposing a symmetry-aware ordering and leveraging the resulting structure of the operators (see \cref{eq:split_s}).
Hence, it is not necessary to build exactly symmetric meshes, and a significant amount of memory and computational resources are saved.

\begin{figure}[ht]
  \centering
  \includegraphics[scale=1]{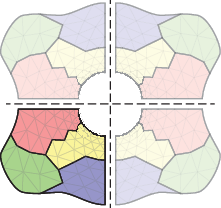}
  \caption{Adequate partitioning of a mesh with 2 reflection symmetries.}\label{fig:dompart}
\end{figure}

Similarly, to replace \xspmv~with \xspmm~effectively, it is required to apply a consistent domain partitioning.
Namely, to distribute the \emph{base mesh} among the available computing resources and extend such a partitioning to the remaining subdomains by the symmetries.
\Cref{fig:dompart} illustrates the above procedure on an arbitrary 2D grid.

The proposed methods have been implemented on top of Chronos~\cite{IsoFriSpiJan21}, a sparse linear algebra library designed for parallel computers.
Chronos provides iterative solvers for linear systems and eigenproblems and advanced preconditioners based on approximate inverses and AMG.
The library is written in C++ with a strongly object-oriented design to ease its development and maintenance.
It presents a hybrid programming model using MPI for inter-node communication and OpenMP and CUDA to take advantage of manycore processors and GPU accelerators, respectively.

As to the preconditioners, Chronos provides both static and adaptive pattern FSAI (sFSAI or aFSAI), which can be used as standalone preconditioners or as smoothers within a multigrid hierarchy.
Chronos also provides an AMG preconditioner implementing a {\em classical} coarsening (that is, a division of the nodes into fine and coarse with no aggregation) and several interpolation schemes.
Namely, classical and ExtI for Poisson-like problems and DPLS for elasticity problems.
Additionally, AMG's quality is improved by using energy minimisation~\cite{JanFraSchOls23} and its application cost reduced through prolongation filtering.

\section{Experimental results}\label{sec:experiments}

This section investigates the advantages of using AMGR to leverage reflection and translational symmetries in industrial applications.
All the \xcfd{} cases considered are governed by the incompressible Navier-Stokes and the continuity equations:
\begin{equation}\label{eq:NS}
	\frac{\partial u}{\partial t} + (u \cdot \nabla) u = \nu \Delta u - \frac{1}{\rho} \nabla p \;\;\;\;\;\;\;
	\nabla \cdot u = 0,
\end{equation}
where $\rho$ and $\nu$ are the density and kinematic viscosity, and $u$ and $p$ are the velocity and pressure fields, respectively.

Regardless of the strategy applied for solving the pressure-velocity coupling of \cref{eq:NS}, a Poisson equation arises, and its solution represents the most computationally intensive part of the simulations.
Considering a classical fractional step projection method~\cite{Chorin1968}, at each time iteration, one computes a predictor velocity, $u^*$, to later project it onto a divergence-free space.
This is done through the gradient of pressure, which is obtained by solving the following Poisson equation:
\begin{equation}\label{eq:CFDpoisson}
  \Delta p= \frac{\rho}{\Delta t}\nabla \cdot u^*,
\end{equation}
where $\Delta t$ is the time-step following the time integration.
See \cite{Trias2014,Verstappen2003} for further details about the discretisation employed.

\begin{figure}[tbp]
	\centering
	\begin{subfigure}{0.35\linewidth}
		\centering
		\includegraphics[trim = 25cm 11cm 28cm 15cm, clip, width=\linewidth]{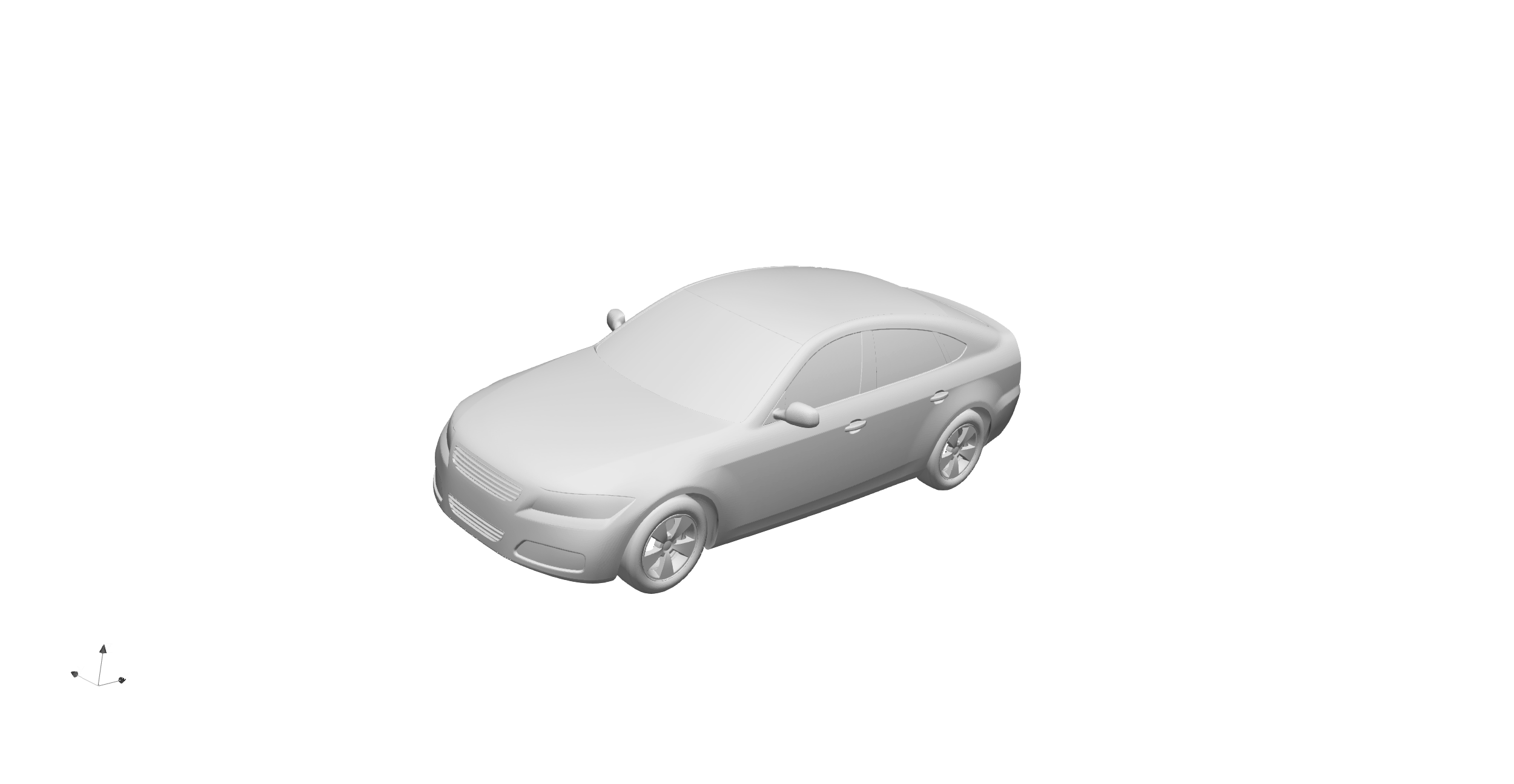}
	\end{subfigure} \qquad
	\centering
	\begin{subfigure}{0.35\linewidth}
		\centering
		\includegraphics[trim = 21cm 11cm 23cm 18cm, clip, width=\linewidth]{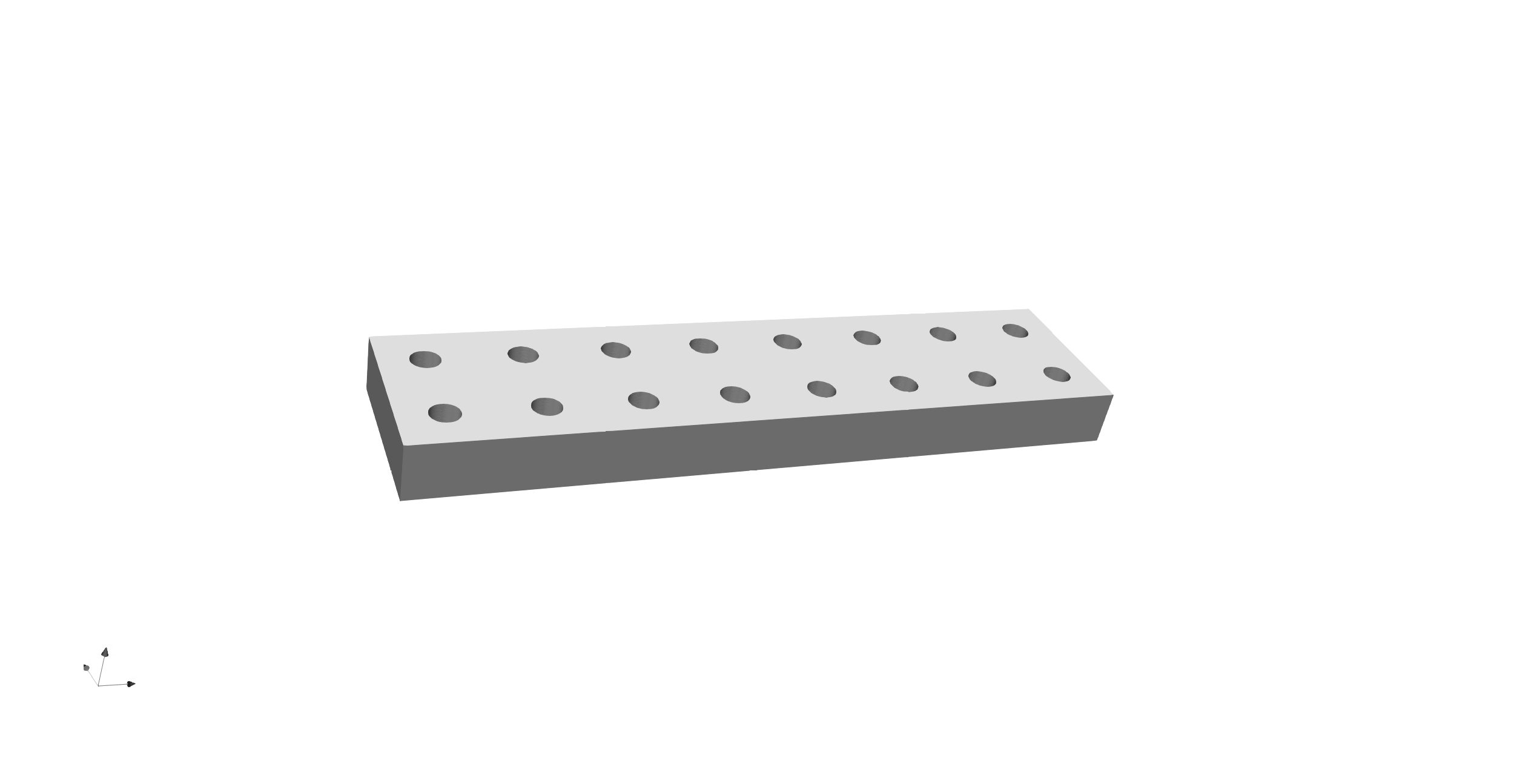}
	\end{subfigure} \\
	\vspace*{0.2cm}
	\begin{subfigure}{0.7\linewidth}
		\centering
		\includegraphics[width=\linewidth]{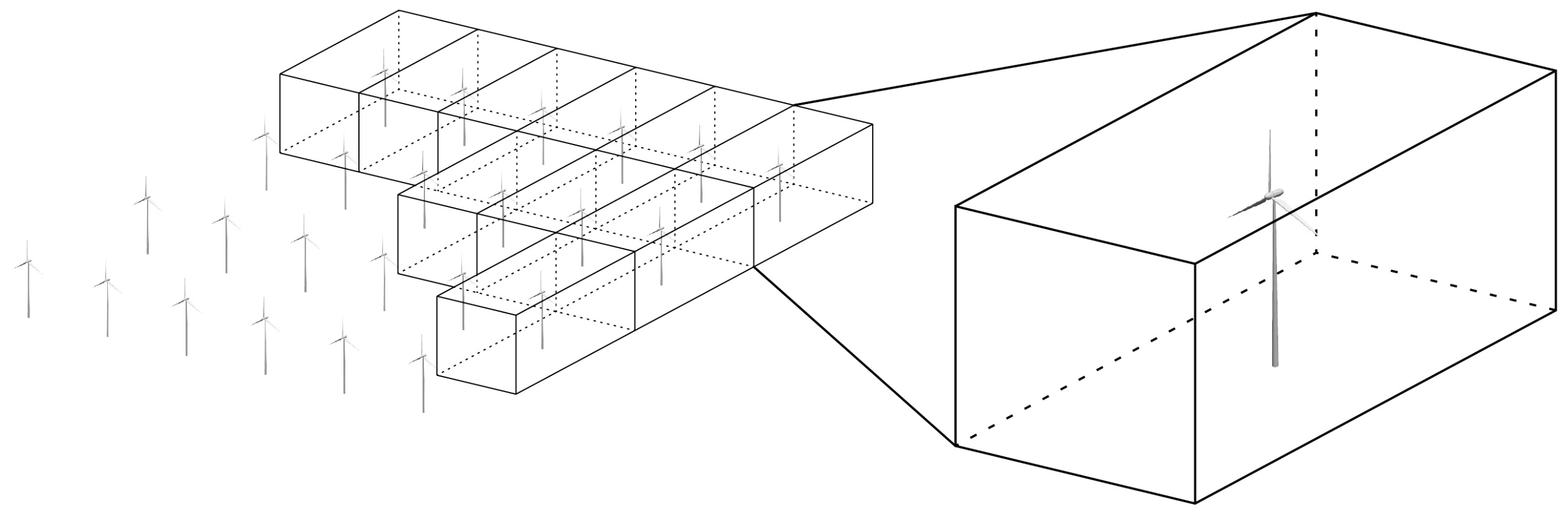}
	\end{subfigure}
	\caption{Industrial test cases used in the numerical experiments. Top: DrivAer car model and finned-tube heat exchanger. Bottom: $6 \times 4$ wind farm.}\label{fig:test_cases}
\end{figure}

The numerical experiments consist of solving \cref{eq:CFDpoisson} using the three industrial test cases of \cref{fig:test_cases}: the flow around a realistic car model, a finned-tube heat exchanger, and the simulation of a wind farm.
While the first two allow for exploiting reflection symmetries, the last allows for studying the further advantages of translational and rotational symmetries.
All the executions rely on combined MPI and multithreaded parallelism and have been conducted on the JFF cluster at the Heat and Mass Transfer Technological Center.
Its non-uniform memory access (NUMA) nodes are equipped with two Intel Xeon 6230 CPUs (20 cores, 2.1 GHz, 27.5 MB L3 cache and 140 GB/s memory bandwidth) linked to 288GB of RAM and interconnected through 7 GB/s FDR Infiniband.

The numerical experiments compare AMGR- to AMG-preconditioned PCG.
For simplicity, the right-hand side used in the tests is a random vector, which ensures a broad spectrum of frequencies in the resulting residual, allowing, in turn, for a complete assessment of the preconditioner's effectiveness in damping all error components.
AMG uses the standard configuration for CFD problems discussed in \cref{sec:impl}.
That is, PMIS coarsening~\cite{Sterck2006} and ExtI interpolation~\cite{DeSFalNolYan08}.
The smoother corresponds to a lightweight aFSAI with five nonzeros per row, and a direct method solves the coarsest level.
Regarding AMGR, as discussed in \cref{ssec:MGR}, the first level of the hierarchy corresponds to the multigrid reduction, and subsequent levels are purely algebraic and use the same configuration as AMG.
For this reason, the coarsening ratio (c-rat) and average nonzeros per row (avg nnzr) in \cref{tab:popul,tab:prolo,tab:car,tab:fintube,tab:windfarm} correspond to the operator at the first level of the multigrid hierarchy, \ie~after the first coarsening, which, in the case of AMGR, is induced by the reduction framework, whereas for AMG by the PMIS coarsening and ExtI interpolation.
Additionally, $n_b$ denotes the number of subdomains in which the domain can be decomposed, \ie~the number of repeated blocks in $\bK = \Id_{n_b}\otimes\K$.
Note that $n_b=1$ corresponds to AMG and $n_b>1$ to AMGR.

As discussed in \cref{ssec:MGR}, the coarsening, interpolation, and smoother used in the multigrid reduction can be adapted to the problem at hand.
In all the CFD problems considered, the smoother will be an aFSAI with five nonzeros per row.
When it comes to the coarsening, \cref{alg:popul} allowed for setting the maximum interpolation distance, $k$.
\Cref{tab:popul} illustrates the impact of $k$ on an AMGR using DPLS with energy minimisation (DPLS + EMIN).
As noted in \cref{ssec:MGR}, the larger the $k$, the more aggressive the coarsening but the potentially more inaccurate the interpolation, leading to faster iterations but slower convergences.
Notably, increasing $k$ reduces the number of coarse nodes, making the prolongation and, consequently, the reduced operator sparser.
\Cref{tab:popul} confirms that $k=1$ retains excessive connections, as the number of iterations is practically the same as for $k=2$.
On the contrary, using $k>2$ quickly degrades convergence, with $k=2$ generally being the best-performing value.

\begin{table}[tbp]
{\footnotesize
  \caption{Impact of the interpolation distance on the heat exchanger problem.}\label{tab:popul}
  \begin{center}
  	\begin{tabularx}{\linewidth}{|c| *{4}{ *{2}{Y|}C{0.45cm}| }}
	  \cline{2-13}
  	   \multicolumn{1}{c|}{} 
	   & \multicolumn{3}{c|}{$k=1$}
	   & \multicolumn{3}{c|}{$k=2$}
	   & \multicolumn{3}{c|}{$k=3$}
	   & \multicolumn{3}{c|}{$k=4$}\\
  	\hline
	   \multirow{2}{*}{$n_b$} &
	   \multirow{2}{*}{c-rat} & avg nnzr & \multirow{2}{*}{its} &
	   \multirow{2}{*}{c-rat} & avg nnzr & \multirow{2}{*}{its} &
	   \multirow{2}{*}{c-rat} & avg nnzr & \multirow{2}{*}{its} &
	   \multirow{2}{*}{c-rat} & avg nnzr & \multirow{2}{*}{its} \\
  	\hline
    1 & 0.36 & 14.5 & 20 & na   & na   & na & na   & na   & na & na   & na   & na  \\
    2 & 0.36 & 50.8 & 18 & 0.14 & 37.2 & 19 & 0.07 & 19.9 & 30 & 0.04 & 12.1 & 907 \\
    4 & 0.36 & 50.6 & 17 & 0.15 & 37.7 & 18 & 0.08 & 21.3 & 27 & 0.05 & 14.7 & 427 \\
    8 & 0.37 & 49.5 & 16 & 0.15 & 37.6 & 18 & 0.09 & 22.6 & 26 & 0.06 & 16.9 & 283 \\
  	\hline
	  \end{tabularx}
  \end{center}
}
\end{table}

The interpolation strategy is critical to make AMGR converge comparably to the underlying \xamg.
Using $k=2$ ensures a roughly twice faster coarsening but requires robust long-distance interpolations to preserve convergence.
\Cref{tab:prolo} compares the multigrid reductions arising from ExtI and DPLS with and without energy minimisation.
Using $k=2$ leads to fine nodes only surrounded by fine nodes, and in such situations, ExtI is not as effective~\cite{DeSFalNolYan08}.
For this reason, in \cref{tab:prolo}, it uses $k=1$.
Unfortunately, neither ExtI produces an accurate interpolation nor does energy minimisation improve it, as the energy of the prolongation is already low, rendering ExtI ineffective for the multigrid reduction.
On the other hand, the dynamic approach of DPLS made it robust enough to use $k=2$.
DPLS was designed to accommodate multiple test vectors, which makes it perform poorly on CFD problems due to their one-dimensional near-null space~\cite{PalFraJan19}.
However, \cref{tab:prolo} demonstrates how effectively energy minimisation improves the quality of DPLS interpolations while keeping the reduced operator sparse, resulting in an AMGR that converges comparably to AMG despite the aggressive coarsening.
Hence, in the remainder of this section, the multigrid reduction will use $k=2$ and DPLS with energy minimisation.

\begin{table}[tbp]
{\footnotesize
  \caption{Impact of the interpolation scheme on the heat exchanger problem.}\label{tab:prolo}
  \begin{center}
  	\begin{tabularx}{\linewidth}{|c|*{4}{Y|C{0.45cm}|C{0.6cm}|}}
	  \cline{2-13}
  	   \multicolumn{1}{c|}{} 
	   & \multicolumn{3}{c|}{ExtI}
	   & \multicolumn{3}{c|}{ExtI + EMIN}
	   & \multicolumn{3}{c|}{DPLS}
	   & \multicolumn{3}{c|}{DPLS + EMIN}\\
  	\hline
	   \multirow{2}{*}{$n_b$} &
	   avg nnzr & \multirow{2}{*}{its} & \multirow{2}{*}{t-sol} &
	   avg nnzr & \multirow{2}{*}{its} & \multirow{2}{*}{t-sol} &
	   avg nnzr & \multirow{2}{*}{its} & \multirow{2}{*}{t-sol} &
	   avg nnzr & \multirow{2}{*}{its} & \multirow{2}{*}{t-sol} \\
  	\hline
    1 & 14.5 & 20  & 1.54 & na    & na  & na   & na   & na  & na   & na   & na & na   \\
    2 & 30.7 & 175 & 8.79 & 111.3 & 140 & 9.63 & 16.3 & 810 & 38.2 & 37.2 & 19 & 1.09 \\
    4 & 30.5 & 77  & 3.42 & 110.7 & 59  & 4.22 & 16.8 & 303 & 14.0 & 37.7 & 18 & 0.95 \\
    8 & 30.0 & 58  & 2.40 & 108.1 & 46  & 3.05 & 19.3 & 192 & 8.27 & 37.6 & 18 & 0.89 \\
  	\hline
	  \end{tabularx}
  \end{center}
}
\end{table}

The top of \cref{fig:test_cases} displays the two test cases used to assess the performance of AMGR when exploiting reflection symmetries.
Namely, the simulation of the flow around the DrivAer fastback car model~\cite{Heft2012} and within a finned-tube heat exchanger, on which we exploit one and three reflection symmetries, respectively.
The results for the DrivAer case in \cref{tab:car} make clear the advantages of AMGR, which preserves the excellent convergence of the standard \xamg~and yields 43\% speed-ups thanks to the combined action of the aggressive coarsening and the faster top-level smoother of \cref{eq:topsmoother}, which leverages \xspmm.
The larger the number of repeated blocks in $\bK = \Id_{n_b}\otimes\K$, the greater the benefits of replacing \xspmv~with \xspmm.
In this sense, the heat exchanger problem extends the results of \cref{tab:car} up to 3 symmetries.
According to \cref{tab:fintube}, regardless of $n_b$, the aggressive coarsening induced by the multigrid reduction results in almost the same coarse operator, $\A_c$.
Indeed, both the coarsening rate and density remain almost constant.
In fact, when increasing $n_b$, the multigrid reduction becomes slightly less aggressive.
Such behaviour follows from the fact that \cref{alg:popul} populates the interface with roughly the same inner unknowns, MIS($T^k$), regardless of $n_b$, whereas the interface grows (slightly) with $n_b$.
However, larger $n_b$ implies faster top-level smoothing, resulting in up to 73\% accelerations.

\begin{table}[htbp]
{\footnotesize
  \caption{DrivAer problem with 106.4M unknonws on five JFF nodes.}\label{tab:car}
  \begin{center}
    \begin{tabular}{|c|c|c|c|c|c|c|} \hline
      preconditioner & $n_b$ & c-rat & avg nnzr & its & t-sol (s) & speed-up \\ \hline
      AMG  & 1 & 0.36 & 14.7 & 26 & 7.71 & 1.00 \\
      AMGR & 2 & 0.14 & 37.4 & 26 & 5.39 & 1.43 \\ \hline
    \end{tabular}
  \end{center}
}
\end{table}

\begin{table}[htbp]
{\footnotesize
  \caption{Heat exchanger problem with 18.4M unknonws on two JFF nodes.}\label{tab:fintube}
  \begin{center}
    \begin{tabular}{|c|c|c|c|c|c|c|} \hline
      preconditioner & $n_b$ & c-rat & avg nnzr & its & t-sol (s) & speed-up \\ \hline
      AMG  & 1 & 0.36 & 14.5 & 20 & 1.54 & 1.00 \\
      AMGR & 2 & 0.14 & 37.5 & 19 & 1.09 & 1.41 \\
      AMGR & 4 & 0.15 & 37.4 & 18 & 0.95 & 1.62 \\
      AMGR & 8 & 0.15 & 37.6 & 18 & 0.89 & 1.73 \\ \hline
    \end{tabular}
  \end{center}
}
\end{table}

The results in \cref{tab:car,tab:fintube} make clear the potential of applying AMGR on domains arising from recurring structures, such as the $6\times4$ wind farm in \cref{fig:test_cases}.
Indeed, while reflection symmetries are generally restricted to $n_b\le 8$, rotational and translational symmetries are not, allowing, for instance, the simulation of an arbitrarily large wind farm by just discretising a single wind turbine.
\Cref{tab:windfarm} summarises the results obtained on a $6\times4$ wind farm discretised according to \cite{Calaf2010}.
The wind turbines are introduced through the immersed boundary method and, therefore, discretised with a structured grid stretched around the blades and coarsened vertically towards the atmospheric boundary layer.
The first thing to note is that the optimal number of subdomains into which dividing the wind farm is not $n_b=24$.
In fact, relatively small values of $n_b$ led to maximum speed-ups, as the application of \xamg~on the reduced operator quickly counterbalanced \xspmm's accelerations.
This follows from the aforementioned balance between the effectiveness of the multigrid reduction and \xspmm{} gains: the larger $n_b$, the faster \xspmm, but the larger $\nifc$.
The wind farm problem is susceptible to such a trade-off, as the growing coarsening ratios in \cref{tab:windfarm} indicate.
Nevertheless, AMGR also proved effective on recurring structures, resulting up to 70\% faster than the underlying \xamg.

The quick saturation of \xspmm{} gains as $n_b$ grows is due to the very high sparsity of the coefficient matrices considered, which makes the best-performing top-level smoother very lightweight, only having five nonzeros per row.
Despite obtaining significant speed-ups, this compromised the advantages of AMGR, not only by making \xspmm's acceleration quickly counterbalanced by the application of \xamg~on the reduced operator but also by making both \xamg~and AMGR have comparable memory footprints.
On top of that, note that the denser the matrices, the more data accesses are saved by \xspmm{}, boosting its accelerations.
In this sense, higher-order schemes or, more generally, applications entailing a denser Poisson's equation (\eg~linear elasticity or geomechanical problems) would strengthen the benefits of \xspmm~and, therefore, of AMGR.

\begin{table}[htbp]
{\footnotesize
  \caption{Wind farm problem with 52M unknonws on two JFF nodes.}\label{tab:windfarm}
  \begin{center}
    \begin{tabular}{|c|c|c|c|c|c|c|} \hline
      preconditioner & $n_b$ & c-rat & avg nnzr & its & t-sol (s) & speed-up \\ \hline
      AMG  & 1  & 0.41 & 13.6 & 37 & 11.0 & 1.00 \\
      AMGR & 2  & 0.23 & 22.1 & 35 & 8.06 & 1.36 \\
      AMGR & 3  & 0.23 & 22.0 & 32 & 6.48 & 1.70 \\
      AMGR & 4  & 0.25 & 20.4 & 36 & 7.54 & 1.46 \\
      AMGR & 6  & 0.24 & 21.7 & 38 & 8.50 & 1.29 \\
      AMGR & 8  & 0.25 & 20.2 & 38 & 8.74 & 1.26 \\
      AMGR & 12 & 0.25 & 20.7 & 37 & 8.78 & 1.25 \\
      AMGR & 24 & 0.26 & 19.9 & 42 & 11.0 & 1.00 \\ \hline
    \end{tabular}
  \end{center}
}
\end{table}

\begin{figure}[tbp]
	\centering
	\begin{subfigure}{0.495\linewidth}
		\centering
    \includegraphics[width=\linewidth]{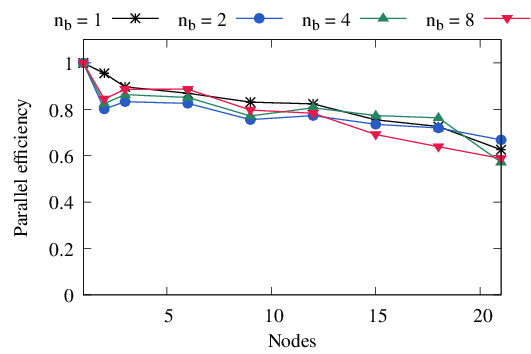}
	\end{subfigure} \hfill
	\centering
	\begin{subfigure}{0.495\linewidth}
  \includegraphics[width=\linewidth]{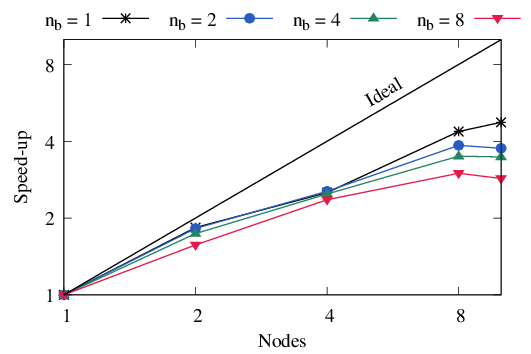}
	\end{subfigure}
	\caption{AMG and AMGR scaling on the unit cube. Left: weak scaling with a workload of $384^3$ unknowns per node. Right: strong scaling on a $256^3$ mesh.}\label{fig:scaling}
\end{figure}

When it comes to the scalability of AMGR, it is worth recalling \cref{sec:impl}.
As discussed, to replace \xspmv~with \xspmm, it is necessary to apply a consistent domain partitioning.
That is, instead of distributing the rows of $\bK = \Id_{n_b}\otimes\K$ among all the available resources, we need to distribute the smaller $\K$ so that its coefficients are effectively reused in the \xspmm~(analogously for the top-level smoother).
Of course, distributing the smaller sub-matrix entails potentially larger communications.
According to the results of \cref{fig:scaling}, this is not a severe issue.
On the one hand, AMGR preserves the excellent weak scalability of \xamg.
On the other, even if AMGR's strong scalability is affected by the extra communication overheads, this effect is far from critical, especially considering the relatively large workloads that extreme-scale simulations typically entail.

\section{Conclusions}\label{sec:concl}

This paper presented a multigrid reduction framework for accelerating \xamg~on regular domains.
Given that modern supercomputers prioritise FLOP performance, most computational physics applications are memory-bound and, therefore, unable to achieve hardware's theoretical peak throughput.
To address this limitation, we first showed that given an arbitrarily complex geometry presenting reflection, translational or rotational symmetries, it is possible to apply a consistent ordering that makes the coefficient matrix (and preconditioners) satisfy a regular block structure.
This, in turn, allows the standard \xspmv~to be replaced with the more compute-intensive \xspmm, accelerating the application of the preconditioners.
This paper focused on accelerating \xamg~owing to its numerical and computational effectiveness.
However, it is worth noting that the intermediate strategies presented naturally apply to other preconditioners. 

AMGR, the proposed multigrid reduction framework, introduces an aggressive coarsening to the first level of the multigrid hierarchy, reducing the memory footprint, setup and application costs of the top-level smoother.
While preserving the excellent convergence of \xamg, replacing \xspmv~with \xspmm~yielded significant speed-ups.
Remarkably enough, AMGR, the resulting preconditioner, does not have any requirement on the boundary conditions, which can be asymmetric.
Hence, not only it succeeds in outperforming the standard \xamg~but also applies to a significantly broader range of industrial applications.
Numerical experiments on industrial \xcfd~applications demonstrated up to 70\% speed-ups in the solution of Poisson's equation with AMGR instead of \xamg.
Furthermore, even if AMGR potentially entails higher communication overheads, strong and weak scalability analyses revealed no serious degradation compared to \xamg.

The coefficient matrices arising from the incompressible \xcfd{} applications considered are remarkably sparse and, as a result, so is the best-performing top-level smoother.
Indeed, best results were obtained smoothing with a very light \xfsai, and even if we obtained significant speed-ups, this compromised the advantages of AMGR.
Firstly, leveraging relatively few symmetries led to maximum speed-ups, as  \xspmm's accelerations were quickly counterbalanced by the application of \xamg~on the reduced operator, whose size grows slightly with the number of subdomains.
Secondly, both \xamg~and AMGR had a comparable memory footprint.
Again, this follows from the particularly lightweight smoother required by the cases considered.
In this sense, higher-order schemes or, more generally, applications entailing a denser Poisson's equation (\eg~linear elasticity or geomechanical problems) would strengthen the advantages of \xspmm~and, therefore, AMGR.

For all these, immediate lines of work include applying AMGR to (denser) problems arising from structural mechanics.
Additionally, we plan to optimise the setup phase, which was irrelevant for the incompressible \xcfd~applications considered, in which the coefficient matrix generally remains constant over time, making relatively large pre-processing overheads acceptable.
Significant accelerations are expected from computing the top-level smoother on the base mesh instead of the entire domain.
Finally, we want to tackle the GPU implementation of AMGR and study the extension of the multigrid reduction framework towards arbitrary domains, which, despite not exploiting \xspmm, would benefit from the aggressive coarsening.


\bibliographystyle{siamplain}
\bibliography{library}
\end{document}